\documentclass[12pt,a4paper,oneside,notitlepage]{amsart}
\usepackage{amssymb}
\usepackage[arrow,matrix,curve]{xy}\SilentMatrices
\def\xyma{\xymatrix@M.7em}

\usepackage[centertags]{amsmath}

\usepackage{epsf,epsfig}
\usepackage{latexsym, amsbsy, amsthm, amscd}

\usepackage{amsfonts}
\usepackage{amscd}
\usepackage{euscript}
\usepackage{amsmath}

\let\dlim\underrightarrow
\numberwithin{equation}{section}
\newtheorem{theorem}{Theorem}[section]

\newtheorem{prop}[theorem]{Proposition}
\newtheorem{lemma}[theorem]{Lemma}
\newtheorem{cor}[theorem]{Corollary}

\begin{document}
\title{Homotopy types of reduced 2-nilpotent simplicial groups}

\author[Baues]{Hans-Joachim Baues}
\address{Max-Planck-Institut fur Mathematik, Vivatsgasse 7, 53111, Bonn, Germany}
\email{baues@mpim-bonn.mpg.de}

\author[Mikhailov]{Roman Mikhailov}
\address{Steklov Mathematical Institute, Gubkina 8, 119991 Moscow, Russia}
\email{romanvm@mi.ras.ru}
\begin{abstract}
We classify the homotopy types of reduced 2-nilpotent simplicial
groups in terms of the homology an
d boundary invariants $b,\beta$.
This contains as special cases results of J.H.C. Whitehead on
1-connected 4-dimensional complexes and of Quillen on reduced
2-nilpotent rational simplicial groups. Moreover it yields for
1-nilpotent (or abelian) simplicial groups a classification due to
Dold-Kan. Our result describes a new natural structure of the
integral homology of any simply connected space. We also classify
the homotopy types of connective spectra in the category of
2-nilpotent simplicial groups. Moreover we compute homotopy groups
of spheres in the category of $m$-nilpotent groups for $m=2,3$ and
partially for $m=4,5$.
\end{abstract}
\maketitle

A classical result of D.M. Kan shows that the homotopy type of a
connected space is determined by a simplicial group, see
\cite{Kan}. E. Curtis \cite{Curtis:63} observed that for a
1-connected $n$-dimensional space the associated simplicial group
can be chosen to be $m$-nilpotent where $n\leq 2+\{log_2(m+1)\}.$
Here $\{n\}$ is the least integer  $\geq n$. Hence a 1-connected
4-dimensional space is equivalent to a reduced 2-nilpotent group
and such spaces were classified by J.H.C. Whitehead
\cite{Whitehead} in terms of the homology groups $H_2,H_3,H_4$ , a
boundary $b: H_4\to \Gamma H_2$ and an invariant $\beta\in
Ext(H_3,coker(b))$ given by the homotopy group $\pi_3$. Moreover,
a result of D. Quillen \cite{QuillenRational} on reduced rational
nilpotent simplicial groups shows that reduced rational
2-nilpotent groups are classified by a boundary $b: B\to [B,B]$ of
degree $-1$ where $B$ is a graded $\mathbb Q$-vector space with
$B_i=0$ for $i\leq 0$. Here $[B,B]$ is defined in the free Lie
algebra $L(B)$ generated by $B$ so that
$$
[B,B]=\bigoplus_{i\ \text{odd}} \Gamma B_i\oplus\bigoplus_{i\
\text{even}}\Lambda^2B_i\oplus \bigoplus_{i<j}B_i\otimes B_j
$$
where $[B_i,B_i]=\Gamma B_i$ (resp. $=\Lambda^2B_i$). In fact, $b$
is part of the differential in a minimal model constructed by
Baues-Lemaire \cite{BL}.

On the other hand, 1-nilpotent (or abelian) simplicial groups are
equivalent to chain complexes by the classical Dold-Kan theorem.
Various authors extended the Dold-Kan theorem for more general
simplicial groups (see \cite{Conduche} and \cite{Carrasco-Cegarra}
and for 2-nilpotent groups the thesis of M. Hartl \cite{Hartl1})
without, however, achieving a description of homotopy types.

Knowing about this background in the literature there is a clear
motivation to study the homotopy types of reduced 2-nilpotent
simplicial groups. In this paper we obtain their complete
description in terms of graded abelian groups and boundary
invariants $b$ and $\beta$. This classification contains as
special cases the results of J.H.C. Whitehead and D. Quillen
mentioned above. For the description of the boundary invariants we
need quadratic functors $\Gamma$ and $\Lambda^2$ and quadratic
torsion functors $\Omega, R$ which were used by Eilenberg-MacLane
\cite{EM} and Cartan \cite{Cartan} for the computation of the
homology of an Eilenberg-MacLane space.

We also consider the homotopy types of connective spectra in the
model category of 2-nilpotent simplicial groups. They are
classified in a particular simple fashion by $\mathcal
F$-algebras, where $\mathcal F$ is the free $\mathbb Z_2$-algebra
generated by elements $Sq_k^{nil}$ for $k$ even $>0$, see section
11.

Since any simplicial group $G$ yields a 2-nilpotent simplicial
group $G/\gamma_3(G)$ by dividing out triple commutators we see
that we can apply our result to $G/\gamma_3(G)$. Therefore we
obtain by the invariants $b,\beta$ a new natural structure of the
integral homology of any 1-connected space. This structure
generalizes the well known action of the Steenrod algebra. In
fact, the action of $Sq_k^{nil}\in \mathcal F$ corresponds to the
dual of the Steenrod square $Sq^k$.

We also describe some aspects of the homotopy theory of
2-nilpotent simplicial groups, in particular, we compute all
homotopy groups of spheres in this category. In an Appendix, we
determine homotopy groups of spheres in the category of
$m$-nilpotent simplicial groups for $m=3$ and partially for
$m=4,5$. \vspace{.5cm}
\section{Simplicial objects and chain complexes}
\vspace{.5cm} Let $\sf C$ be a category with an object $0$ which
is initial and final in $\sf C$. For objects $X,Y$ in $\mathcal
C$, one has the unique zero morphism $0: X\to 0\to Y$. A {\it
graded object} $X_*$ in $\sf C$ is a sequence of objects
$X_*=\{X_n,\ n\in \mathbb Z\}$ in $\sf C$. A {\it chain complex}
$(X_*,d)$ is a graded object together with boundary morphisms $d:
X_n\to X_{n-1}$ satisfying $dd=0$. A {\it simplicial object} in
$\sf C$ is a sequence $\{X_n,\ n\geq 0\}$ together with morphisms
$d_i: X_n\to X_{n-1}, s_i: X_n\to X_{n+1},\ i=0,\dots, n$,
satisfying the usual simplicial identities. Such objects are
$r-reduced$ if $X_n=0$ for $n<r$. Let
$$\sf C_*,\ \ Chain(C),\ \ s C$$ be the
categories of graded objects, chain complexes and simplicial
objects respectively. Moreover, let
\begin{equation}
{\sf C}_r,\ \ {\sf Chain(C)}_r,\ \ ({\sf sC})_r
\end{equation}
be the corresponding full subcategories consisting of $r$-reduced
objects. Dold-Kan theorem shows that for an abelian category $\sf
A$ (for example ${\sf A}=Mod(R)$, the category of $R$-modules for
a ring $R$), one has the equivalence of categories
$$
N: ({\sf sA})_r\buildrel{\sim}\over{\longrightarrow} \sf
{Chain(A)}_r,\ r\geq 0.
$$
Here $N$ is the {\it normalization functor} which maps the
simplicial object $A$ in $\sf A$ to the Moore chain complex $N(A)$
with
\begin{equation}\label{ref1.2}
N_q(A)=\begin{cases} \bigcap_{i>0}ker\{d_i: A_q\to A_{q-1}\},\ \text{if}\ q>0,\\
A_0,\ \text{if}\ q=0.
\end{cases}
\end{equation}
The boundary $d: N_q(A)\to N_{q-1}(A)$ is induced by $d_0$ with
$d_0=0$ for $q=0$. Let $N^{-1}$ be the {\it inverse} of the
normalization functor. We define the {\it homotopy group} by the
quotient in $\sf A$:
\begin{equation}\label{ref1.3}
\pi_q(A)=\frac{ker\{d:N_qA\to N_{q-1}A\}}{im\{d:N_{q+1}A\to
N_qA\}}
\end{equation}
which is the $q$-th homology $H_qN(A)$ of the chain complex
$N(A)$. The equivalence $N$ shows that a simplicial object in $\sf
A$ can be identified with a (non-negative) chain complex in $\sf
A$.

Let ${\sf Gr}$ be the {\it category of groups} then one has a
normalization functor
$$
N: ({\sf sGr})_r\to {\sf Chain(Gr)}_r,\ r\geq 0.
$$
The functor $N$ carries the simplicial group $A$ to the chain
complex $N(A)$ defined in same way as in (\ref{ref1.2}). Moreover
we can define the (Moore) homotopy group $\pi_q(A)$ of the
simplicial group $A$ by the quotient group in (\ref{ref1.3}) since
$im(d)$ is normal in $ker(d)$.

{\it Weak equivalences} in the categories $\sf Chain(A)$, resp.
$\sf sA$ and $\sf sGr$, are maps which induce isomorphisms for the
homology functor $H_*$, resp. for the homotopy group functor
$\pi_*$. Given a category $\sf C$ with weak equivalences we obtain
the {\it homotopy category} $\sf {Ho (C)}$ which is the
localization of $\sf C$ with respect to the class of weak
equivalences. If $\sf C$ is a model category or a {\it cofibration
category}, the homotopy category $\sf Ho(C)$ is well-defined, see
Quillen \cite{HA} and Baues \cite{AH}. For objects $X,Y$ in $\sf
C$, let
\begin{equation}
[X,Y]=[X,Y]_{\sf C}
\end{equation}
be the set of morphisms $X\to Y$ in the homotopy category $\sf
Ho(C)$. A {\it homotopy type} in $\sf C$ is the equivalence class
of an object $X$ in the homotopy category $\sf Ho(C)$.
\vspace{.5cm}
\section{Abelianization and nilization}
\vspace{.5cm} For a group $G$ one has the {\it lower central
series}
$$\dots \subseteq \gamma_{n+1}(G)\subseteq \gamma_n(G)\subseteq
\dots\subseteq \gamma_2(G)\subseteq G,$$ where $\gamma_2(G)$ is
the commutator subgroup. Then the quotient
$$
ab(G)=G/\gamma_2(G)
$$
is the {\it abelianization of} $G$. A group $G$ has {\it
nilpotency degree} $\leq k$ if $\gamma_{k+1}(G)$ is trivial, that
is if all $(k+1)$-fold commutators in $G$ are trivial. In this
case we also call $G$ a {\it nil(k)-group}. We only deal with
$nil(2)$-groups though various concepts below have an obvious
analogue for $nil(k)$-groups, $k\geq 2$. Compare the Appendix. We
call the quotient
$$
nil(G)=G/\gamma_3(G)
$$
the {\it nilization of} $G$. Hence we have functors
$$
{\sf Gr}\buildrel{nil}\over{\to}{\sf Nil}\buildrel{ab}\over{\to}
{\sf Ab},
$$
where $\sf Ab$ is the category of abelian groups and where $\sf
Nil$ is the full subcategory of $\sf Gr$ consisting of
$nil(2)$-groups. Clearly one gets $ab(nil(G))=ab(G)$. If
$G=\langle E\rangle$ is a free group generated by the set $E$,
then $A=ab(G)=\mathbb Z[E]$ is the free abelian group and
$nil(G)=\langle E\rangle_{nil}$ is the free $nil(2)$-group. Here
$\mathbb Z[E]$ and $\langle E\rangle_{nil}$ are free objects in
the categories $\sf Ab$ and $\sf Nil$ respectively. If $G=\langle
E\rangle_{nil}$ is a free $nil(2)$-group we have the central
extension of groups
\begin{equation}\label{seqlambda}
\Lambda^2(A)\buildrel{w}\over{\rightarrowtail} G
\buildrel{p}\over{\twoheadrightarrow} A,
\end{equation}
where $A$ is the abelianization of $G$ and where $w$ is the
commutator map with $w(px\wedge py)=[x,y]:=x^{-1}y^{-1}xy$ for
$x,y\in G$. For this recall that $\Lambda^2(A)$ is the exterior
square of the abelian group given by the quotient
$\Lambda^2(A)=A\otimes A/\{a\otimes a\sim 0\}.$ We point out that
the homomorphisms $w$ and $p$ are natural for $\phi: G\to G'\in
\sf Gr$ that is $p\phi=\phi_*p,\ w\Lambda^2(\phi_*)=\phi w$ with
$\phi_*=ab(\phi)$.

The categories $\sf Gr$ and $\sf Nil$ are closed under limits and
colimits. In fact, colimits in $\sf Nil$ are obtained by
nilization of the corresponding colimits in $\sf Gr$. For example,
the {\it sum} $A\vee B$ in $\sf Nil$ is given by $A\vee
B=nil(A\star B)$, where $A\star B$ is the free product in $\sf
Gr$.

It is useful to note that for any $A,B$ in $\sf Nil$ one has the
functorial short exact sequence (see for example 7.10 in
\cite{QECG})
$$
0\to ab(A)\otimes ab(B) \buildrel{w}\over{\to} A\vee B\to A\times
B\to 0,
$$
where $w(\bar a\otimes \bar b)=[a,b]$. Here $\bar a$ denotes the
class of $a\in A$ in $ab(A)$.

Let $r\geq 0$. Then the functors $nil$ and $ab$ above induce
functors between categories of simplicial groups
$$
({\sf sGr})_r\buildrel{nil}\over{\to} ({\sf
sNil})_r\buildrel{ab}\over{\to} ({\sf sAb})_r.
$$
For $r\geq 1$ these functors carry weak equivalences to weak
equivalences so that the induced functors between homotopy
categories
$$
{\sf Ho(sGr)}_r\buildrel{nil}\over{\to} {\sf
Ho(sNil)}_r\buildrel{ab}\over{\to} {\sf Ho(sAb)}_r.
$$
are well-defined. Using results of Quillen \cite{HA}, we see that
${\sf (sGr)}_r,\ {\sf (sNil)}_r$ and ${(\sf sAb)}_r$ are actually
model categories. Moreover by the Dold-Kan theorem it is well
known that the normalization
$$
N: {\sf Ho(sAb)}_r\buildrel{\sim}\over{\longrightarrow} {\sf
Ho(Chain(Ab))}_r
$$
is an equivalence of homotopy categories. Homotopy types in ${\sf
Ho(Chain(Ab))}_r$ are identified via homology with isomorphism
types of graded abelian groups in ${\sf Ab}_r$. Hence also
homotopy types of simplicial abelian groups in ${\sf s(Ab)}_r$ are
given by graded abelian groups in ${\sf Ab}_r$.

Let $\sf CW$ be the category of CW-complexes $X$ with trivial
$0$-skeleton $X^0=*$. Morphisms are base point preserving maps.
Then homotopy $\simeq$ of such maps yields the quotient category
${\sf CW}/\simeq$ which is the usual homotopy category of
algebraic topology. For $X,Y\in \sf CW$, let $[X,Y]$ be the set of
homotopy classes $X\to Y$ in ${\sf CW}/\simeq$. Let ${\sf CW}_r$
be the full subcategory of $\sf CW$ consisting of CW-complexes $X$
with trivial $(r-1)$-skeleton. Kan \cite{HT} constructed a functor
$G: {\sf CW}_{r+1}\to {\sf (sGr)}_r$ which induces an equivalence
of homotopy categories
$$
G: {\sf CW}_{r+1}/\simeq\ \buildrel{\sim}\over{\longrightarrow}
{\sf Ho(sGr)}_r,\ r\geq 0.
$$
This functor carries a CW-complex $X$ to the {\it Kan loop group}
$G(X)$ which is a free simplicial group (see also Curtis
\cite{SH}). Let
$$
G^{ab}(X)=ab(G(X)),\ \text{resp.}\ G^{nil}(X)=nil(G(X))
$$
be the abelianization, resp. nilization of the Kan loop group
which we also call the {\it abelianization}, resp. the {\it
nilization} of the space $X$. Let $C_*X$ be the {\it cellular
chain complex} of $X$ given by $C_nX=H_n(X^n,X^{n-1})$ and let
$\tilde C_*X=C_*X/C_*(*)$ be the reduced cellular chain complex.
Then one has a weak equivalence
$$
s^{-1}\tilde C_*(X)\buildrel{\sim}\over{\longrightarrow}
NG^{ab}(X)
$$
which is a natural isomorphism in ${\sf Ho(Chain(Ab))}_r$ for
$X\in {\sf CW}_{r+1}$. Here $s^kC$ is the {\it k-fold suspension}
of the chain complex $C$ with $(s^kC)_n=C_{n-k}$ and
$d(s^kx)=(-1)^{k|x|}dx$ for $x\in C,\ k\in \mathbb Z$. The
equivalence shows that the homotopy type of the abelianization of
$X$ coincides with the homology $H_*(X)$ of the space $X$.
\vspace{.5cm}
\section{Derived functors of the exterior square} \vspace{.5cm}
For abelian groups $A,B$ let $A\otimes B$ and $A*B$ be the tensor
product and the torsion product respectively. The tensor product
of abelian groups leads to the notion of the {\it tensor product}
$A\otimes B$ of graded abelian groups $A,B$ with
$$
(A\otimes B)_n=\bigoplus_{i+j=n}A_i\otimes B_j.
$$
We also need the {\it ordered tensor product}
$A\buildrel{>}\over\otimes B$ defined by
$$
(A\buildrel{>}\over{\otimes} B)_n=\bigoplus_{i+j=n,\
i>j}A_i\otimes B_j.
$$
Analogically define the {\it ordered torsion product}
$A\buildrel{>}\over{*}B$ as
$$
(A\buildrel{>}\over{*}B)_n=\bigoplus_{i+j=n,\ i>j}A_i* B_j.
$$

The tensor product, torsion product and the ordered tensor and
torsion product are in the obvious way bifunctors.

Next we use the exterior square $\Lambda^2$ and Whitehead's
quadratic functor $\Gamma$ which are functors from abelian groups
to abelian groups. They define the {\it weak square functor}
$$
sq^\otimes: {\sf Ab_r\to Ab_r}
$$
which is given by
$$
sq^\otimes(A)_n=\begin{cases} \Gamma(A_m),\ \text{if}\ n=2m,\ m\ \text{odd},\\
\Lambda^2(A_m),\ \text{if}\ n=2m,\ m\ \text{even},\\
0,\ \text{otherwise}.\end{cases}
$$
Let $(\mathbb Z_2)_{odd}$ be the graded abelian group which is
$\mathbb Z_2$ in odd degree $\geq 1$ and which is trivial
otherwise, hence $(\mathbb Z_2)_{odd}$ is the reduced homology of
the classifying space $\mathbb RP_\infty=K(\mathbb Z_2,1)$. We now
define the {\it square functor}
\begin{align*}
& Sq^\otimes: {\sf Ab}_r\to {\sf Ab}_r\\
& Sq^\otimes(A)=A\buildrel{>}\over{\otimes}(A\oplus (\mathbb
Z_2)_{odd})\oplus sq^\otimes(A).
\end{align*}
Clearly the square functor is quadratic. The cross-effect is
$$
Sq^\otimes(A|B)=A\otimes B
$$
and one has the operators
$$
Sq^\otimes(A)\buildrel{H}\over{\rightarrow} A\otimes
A\buildrel{P}\over{\rightarrow} Sq^\otimes(A)
$$
which are induced by the diagonal and the folding map
respectively.

Define also the {\it torsion square functor} $$Sq^\star(A):{\sf
Ab}_r\to {\sf Ab}_r$$ by setting
$$
Sq^\star(A)=(A\buildrel{>}\over{*}(A\oplus (\mathbb
Z_2)_{odd}))\oplus sq^\star(A),
$$
where
$$
sq^\star(A)_n=
\begin{cases} \Omega(\pi_mX),\ n=2m,\ m\ \text{even}\\
R(\pi_mX),\ n=2m,\ m\ \text{odd}\\
0,\ \text{otherwise}
\end{cases}
$$
Here $R(A)=H_5K(A,2)$ and $\Omega(A)=H_7K(A,3)/(\mathbb Z_3\otimes
A)$ are functors of Eilenberg-MacLane with
$R(A|B)=\Omega(A|B)=A*B$ and $$R(\mathbb Z_n)=\mathbb Z_{(2,n)},\
\Omega(\mathbb Z_n)=\mathbb Z_n,\ R(\mathbb Z)=\Omega(\mathbb
Z)=0.$$

For developing homotopy theory of 2-nilpotent simplicial groups we
need the description of the derived functors of the exterior
square. They were described in \cite{BauesPirashvili} as follows.
Let $X$ be a simplicial group which is free abelian in each
degree. Then there exists the following short exact sequence of
graded abelian groups
\begin{equation}\label{deri}
0\to Sq^\otimes(\pi_*(X))\to \pi_*(\Lambda^2 X)\to
Sq^\star(\pi_*(X))[-1] \to 0
\end{equation}
where $\pi_*(X)$ and $\pi_*(\Lambda^2 X)$ are the graded homotopy
groups of $X$ and $\Lambda^2 X$ respectively. \vspace{.5cm}
\section{Homotopy theory of 2-nilpotent simplicial groups}
\vspace{.5cm} Homotopy theory in the category $\sf sNil$ of
simplicial $nil(2)$-groups relies on the following definition of
weak equivalences, fibrations and cofibrations. Weak equivalences
are maps inducing isomorphisms on homotopy groups, {\it
fibrations} are the maps $f$ for which $N_qf$ is surjective for
$q>0$ and {\it cofibrations} are retracts of free maps. Here an
injective map $f: X\to Z$ in $\sf sNil$ is {\it free} if there are
subsets $C_q\subset Z_q$ for each $q$ such that $(i)\
\eta^*C_q\subset C_p$ whenever $\eta: [q]\to [p]$ is a surjective
monotone map, $(ii)\ f_q+g_q: X_q\vee FC_q\to Z_q$ is an
isomorphism for all $q$. Here $FC_q=nil\langle C_q\rangle$ is the
free $nil(2)$-group generated by $C_q$ and $g_q: FC_q\to Z_q$ is
the extension of $C_q\subset Z_q$. Let $*$ be the initial object
in $\sf sNil$. We say that $G$ in $\sf sNil$ is free if $*\to G$
is free. Let ${\sf (sNil)}_c$ be the full subcategory of the free
objects in $\sf sNil$. For example, for a space $X$, the
nilization $G^{nil}(X)$ of Kan's loop group is a free object.

\begin{prop}
With these definitions the category $\sf sNil$ of simplicial
$nil(2)$-groups is a closed model category. All objects are
fibrant and free objects form a sufficiently large class of
cofibrant objects. Hence one has a notion of homotopy $\simeq$ in
${\sf (sNil)}_c$ such that the inclusion ${\sf (sNil)}_c\subset
\sf sNil$ induces an equivalence of categories
$$
{\sf (sNil)}_c/\simeq\ \buildrel{\sim}\over{\longrightarrow} \sf
Ho(sNil).
$$
\end{prop}
\begin{proof}
This is a consequence of II.§4 in Quillen \cite{HA}, compare in
particular Remark 4. Moreover we use II.§3 in Baues \cite{AH}.
\end{proof}

The proposition shows that the basic definitions and constructions
of homotopy theory are available in the category $\sf sNil$, for
example we have cylinder objects, suspensions, mapping cones,
cofiber sequences, spectral sequences etc as described in a
cofibration category in Baues \cite{AH}.

\begin{prop}
The nilization functor
$$
nil: \sf sGr \to sNil
$$
carries homotopy push outs to homotopy push outs, that is, $nil$
is a model functor in the sense of Baues \cite{AH} (I.1.10).
\end{prop}

The functor $\pi_n: \sf Ho(sNil)\to Gr$ is a representable functor
in the sense that there is a free object $S(n)\in \sf sNil$ for
$n\geq 0$ and a class $i_n\in \pi_nS(n)$ so that for all $X\in \sf
sNil$, the map
$$
[S(n),X]\to \pi_nX,
$$
given by $f\mapsto \pi_n(f)(i_n)$ is an isomorphism for $n\geq 0$.
This is analogous to the situation for topological spaces; indeed,
these isomorphisms virtually demand that we refer to $S(n)$ as the
$n-sphere$ in $\sf sNil$. We can choose $S(n)$ to be the free
object generated by a single element in degree $n$. We have a
homotopy equivalence
$$
S(n)=G^{nil}(S^{n+1})
$$
so that $S(n)$ is the nilization of the standard $(n+1)$-sphere
$S^{n+1}$. {\it Homotopy groups of spheres}
$$
\pi_{n+k}S(n)=[S(n+k),S(n)]
$$
in the category $\sf Ho(sNil)$ can be computed completely.

\begin{prop}
There are generators $i_n,\ \eta_n,\ \eta_n^k$ such that
$$
\pi_{n+k}S(n)=\begin{cases} \mathbb Zi_n,\ k=0\\
\mathbb Z\eta_n,\ n=k\ \text{odd}\\
\mathbb Z_2\eta_n^k,\ 0<k<n,\ k\ \text{odd}\\
0,\ \text{otherwise}.
\end{cases}
$$
\end{prop}

The computation follows from the description of the spaces
$\Lambda^2 G^{ab}(S^n),$ given by Curtis and Schlezinger
\cite{Sch}, \cite{Curtis:63}: For any $n\geq 1,$ $q\geq 0,$ one
has
$$
\pi_{n+q}\Lambda^2 G^{ab}(S^{n+1})=
\begin{cases}\mathbb Z,\ \text{if}\ n\ \text{is odd
and}\ q=n\\  \mathbb Z_2,\ \text{if}\ q=1,3,\dots, 2[n/2]-1,\\
0\ \text{otherwise}.\\
\end{cases}
$$
This is a part of the more general Theorem \ref{mooreh}, and
easily can be proved using the exact sequence (\ref{deri}).

We call $\eta_n: S(2n)\to S(n)$ with $n$ odd a (generalized) {\it
Hopf map}. In fact, $\eta_1$ is the nilization of the classical
Hopf map $S^3\to S^2$. The iterated suspensions of the Hopf maps
yield the elements $\eta_n^k$, that is,
$$
\eta_n^k=\Sigma^{n-k}\eta_k.
$$
Clearly, the identity $i_n:S(n)\to S(n)$ satisfies $\Sigma
i_n=i_{n+1}$.

For $n,m\geq 1$, one has a map
$$
w: S(n+m)\to S(n)\vee S(m)
$$
which in fact is the nilization of the classical Whitehead product
map $S^{n+m+1}\to S^{n+1}\vee S^{m+1}$. For $X$ in $\sf sNil$ we
thus obtain the {\it Whitehead product}
$$
\pi_n(X)\times \pi_m(X)\to \pi_{n+m}(X),
$$
given by $[x,y]=w^*(x,y)$.

\begin{prop}
The Whitehead product in $\sf sNil$ is bilinear and satisfies
\begin{align*}
& [x,y]=(-1)^{|x||y|+1}[y,x],\\
& [x,x]=\begin{cases} 0\ \text{if}\ |x|\ \text{is even},\\
2\eta_n^*(x)\ \text{if}\ |x|=n\ \text{is odd},
\end{cases}\\
& [x,y]=\eta_n^*(x+y)-\eta_n^*(x)-\eta_n^*(y)\ \text{and}\\
& \eta_n^*(-x)=\eta_n^*(x)\ \text{for}\ |x|=|y|=n\ \text{odd}.
\end{align*}
Moreover, all triple Whitehead products in $\sf sNil$ are trivial.
\end{prop}
\begin{cor}
For a space $X$ in $\sf CW$, the nilization $$nil:
\pi_{n+1}X=[S^{n+1},X]\to [S(n),G^{nil}(X)]=\pi_nG^{nil}(X)$$
carries Whitehead products to Whitehead products. This implies
that the nilization
$$
nil: \pi_{2n+1}(S^{n+1})\to \pi_{2n}S(n)=\mathbb Z\eta_n,\ n\
\text{odd},
$$
coincides with the classical Hopf invariant. Hence, for $n=1,3,7$,
the Hopf map $\eta_n$ in $\sf sNil$ is the nilization of the
classical Hopf maps $S^{2n+1}\to S^{n+1}$ and in this case
$\eta_n^*$ is the nilization of the corresponding suspended Hopf
maps. For $n\neq 1,3,7,$ the elements $\eta_n$ and $\eta_n^k$ are
not in the image of $nil$, however, the element
$2\eta_n=[i_{n+1},i_{n+1}]$, with $n$ odd, is always in the image
of $nil$.
\end{cor}

Recall that a function $f: A\to B$ between abelian groups is
quadratic if $f(a)=f(-a)$ and $[a,b]_f=f(a+b)-f(a)-f(b)$ is
bilinear for $a,b\in A$. Let $\gamma: A\to \Gamma A$ be the
universal quadratic function which defines Whitehead's quadratic
functor $\Gamma: \sf Ab\to Ab$. Then $f$ defines a unique
homomorphism $f^{\Box}:\Gamma A\to B$ with $f^{\Box}\gamma=f$. The
theorem above shows that $\eta_n^*$ is quadratic, hence we obtain
natural transformations $(X\in \sf sNil)$
\begin{align*}
& \Gamma(\pi_nX)\to \pi_{2n}X,\ n\ \text{odd},\\
& \Lambda^2(\pi_nX)\to \pi_{2n}X,\ n\ \text{even}\geq 2,
\end{align*}
which are induced by $\eta_n^*$ and the Whitehead square
respectively. \vspace{.5cm}
\section{Homology and Moore objects in $\sf sNil$}
\vspace{.5cm} Let $X$ be an object in $\sf sNil$. We can choose a
weak equivalence $\bar X\to X,$ where $\bar X$ is free in $\sf
sNil$, that is $*\rightarrowtail \bar
X\buildrel{\sim}\over{\longrightarrow} X$ is a cofibrant model of
$X$. Now we define the {\it chain complex} of $X$ by
$$
C_*X=N(ab \bar X)
$$
and we define {\it homology} and {\it cohomology} of $X$ by this
chain complex, that is, for $n\in\mathbb Z$,
\begin{align*}
& H_n(X,A)=H_nC_*(X)\otimes A,\\
& H^n(X,A)=H^nHom(C_*X,A),
\end{align*}
where $A$ is an abelian group of coefficients. We also need the
{\it pseudo-homology}
$$
H_n(A,X)=[C(A,n),C_*X],
$$
which is a set of homotopy classes of chain maps. Here $C(A,n)$ is
a chain complex of free abelian groups with $H_nC(A,n)=A$ and
$H_jC(A,n)=0,$ for $j\neq n$, for example, $C(A,n)$ is given by a
short free resolution of $A$
$$
0\to C'\buildrel{d_A}\over\rightarrow C\to A\to 0
$$
with $C$ in degree $n$ and $C'$ in degree $n+1$.

For a space $X$ in $\sf CW$, we have canonical natural
isomorphisms
\begin{align*}
& \tilde H_{n+1}(X,A)=H_n(G^{nil}(X),A),\\
& \tilde H^{n+1}(X,A)=H^n(G^{nil}(X),A).
\end{align*}
A {\it Moore object} $M(A,n)_{nil}$ in $\sf sNil$ is a free object
with a single non-vanishing homology group $A$ in degree $n\geq
1$. We clearly have $\Sigma M(A,n)_{nil}=M(A,n+1)_{nil}$ and
$C_*M(A,n)_{nil}=C(A,n)$. The nilization of the Moore space
$M(A,n+1)_{\sf CW}$ in $\sf CW$ can be chosen to coincide with
$M(A,n)_{nil}$. For $A=\mathbb Z$ the Moore object
$M(A,n)_{nil}=S(n)$ is the sphere in $\sf sNil$.

Homotopy groups of Moore objects in $\sf sNil$ are completely
computed in the following result:
\begin{theorem}\label{mooreh}
Let $n\geq 1$, $k\in \mathbb Z$. Then the homotopy groups of Moore
objects in $\sf sNil$ are
$$
\pi_{n+k}M(A,n)_{nil}\buildrel{\Theta}\over{=}\begin{cases} A,\ k=0\\
A\otimes \mathbb Z_2,\ 0<k<n,\ k\ \text{odd}\\
A*\mathbb Z_2,\ 0<k<n,\ k\ \text{even}\\
\Gamma A,\ k=n\ \text{odd}\\
\Lambda^2(A)\oplus A*\mathbb Z_2,\ k=n\ \text{even}\\
R(A),\ k=n+1\ \text{even}\\
\Omega(A),\ k=n+1\ \text{odd}\\
0,\ \text{otherwise}.
\end{cases}
$$
\end{theorem}
\begin{proof}
Consider the short exact sequence of simplicial groups:
\begin{equation}\label{lz} 0\to \Lambda^2abM(A,n)_{nil}\to M(A,n)_{nil}\to
abM(A,n)_{nil}\to 0
\end{equation}
All the cases of the theorem besides the case $k=n\ \text{even}$
follow from (\ref{lz}) and the universal coefficient theorem
(\ref{deri}). For the case $k=n,$ (\ref{deri}) implies the short
exact sequence
$$
0\to \Lambda^2(A)\to \pi_{2n}M(A,n)_{nil}\to A*\mathbb Z_2\to 0
$$
which  a priori splits unnaturally. The natural splitting,
however, follows from the periodicity principle in (\ref{deri})
and the following diagram
$$\xyma{ H_6K(A,3)\ar@{->}[rr]^{d^1} \ar@{=}[dr] & &
\pi_{4}M(A,2)_{nil}\\
& \Lambda^2(A)\oplus A*\mathbb Z_2\ar@{->}[ur] & \\}
$$
where $d^1$ is the boundary map, which is a natural isomorphism,
what follows, for example, from the Curtis spectral sequence
argument (see proof of Theorem \ref{lamq}).
\end{proof}
For each homomorphism $\phi: A\to B$ between abelian groups we can
choose a map $\bar \phi\in [M(A,n)_{nil},M(B,n)_{nil}]$ in $\sf
Ho(sNil)$, which induces $\phi$. The map $\bar\phi$ is not unique.
The induced homomorphism $\pi_{n+k}(\bar \phi),$ however, depends
only on $\phi$ so that $A\mapsto \pi_{n+k}M(A,n)_{nil}$ yields a
functor $\sf Ab\to Ab$.

The isomorphism $\Theta$ is natural with respect to maps
$M(A,n)_{nil}\to M(A',n)_{nil}$. For $k=0$, the isomorphism
$\Theta$ is the Hurewicz isomorphism and for $0<k<n$, $k$ odd, we
have $\Theta(a\otimes 1)=a\eta_n^k$, and for $k=n$.

For an object $X$ in $\sf sNil$, we define {\it homotopy groups
with coefficients in} $A$ by
$$
\pi_n(A,X)=[M(A,n)_{nil},X]_{nil}.
$$
As in topology, one has the {\it universal coefficient sequence}
$$
Ext(A,\pi_{n+1}X)\buildrel\Delta\over\rightarrowtail\pi_{n}(A,X)\buildrel{\mu}\over{\longrightarrow}
Hom(A,\pi_nX),
$$
where $\mu(\alpha)$ is the composition
$A=\pi_nM(A,n)_{nil}\buildrel{\alpha_*}\over{\longrightarrow}\pi_nX$.
\begin{prop}
For Moore objects in $\sf CW$ and $\sf sNil$, the nilization
yields a bijection, $n\geq 1$,
$$
[M(A,n+1)_{\sf CW},M(B,n+1)_{\sf CW}]_{\sf
CW}=[M(A,n)_{nil},M(B,n)_{nil}]_{nil}.
$$
Hence, the homotopy category of Moore objects of degree $(n+1)$ in
${\sf CW}/\simeq$ is equivalent to the homotopy category of Moore
objects of degree $n$ in $\sf Ho(sNil)$.
\end{prop}
\begin{proof}
Both sides are part of universal coefficient sequences which are
isomorphic since via nilization $\pi_{n+k+1}M(B,n+1)_{\sf CW}$ in
$\sf CW$ is the same as $\pi_{n+k}M(B,n)_{nil}$ in $\sf Ho(sNil)$
for $k=0,1.$
\end{proof}
\vspace{.5cm}
\section{Quadratic functors}
\vspace{.5cm} Let $F: {\sf Ab\to Ab}$ be a functor. The {\it
cross-effect} of $F$ is the bifunctor defined as
$$
F(X|Y)=ker\{F(r_1,r_2):F(X\oplus Y)\to F(X)\oplus F(Y)\},\ X,Y\in
{\sf Ab}
$$
where the map $F(r_1,r_2)$ is induced by natural retractions
$r_1:X\oplus Y\to X,\ r_2: X\oplus Y\to Y$. The functor $F$ is
linear if $F(0)=0$ and $F(X|Y)=0$ for all $X,Y\in {\sf Ab}$. The
functor $F$ is {\it quadratic} if $F(0)=0$ and the cross-effect
$F(X|Y)$ is linear in each variable $X$ and $Y$.

A quadratic $\mathbb Z$-module is a diagram of abelian groups
$$
M=(M_e\buildrel{H}\over\rightarrow M_{ee}
\buildrel{P}\over\rightarrow M_e)$$ satisfying $HPH=2H$ and
$PHP=2P$.

Let $M=(M_e\buildrel{H}\over\rightarrow
M_{ee}\buildrel{P}\over\rightarrow M_e)$ be a quadratic $\mathbb
Z$-module.,Then $M$ induces a quadratic functor $A\mapsto A\otimes
M,\ A\in {\sf Ab}$ defined as follows. Given an abelian group $A$,
the abelian group $A\otimes M$ has generators $a\otimes m,\
[a,b]\otimes n,\ a,b\in A,\ m\in M_e,\ n\in M_{ee}$ and relations
\begin{align*}
& (a+b)\otimes m=a\otimes m+b\otimes m+[a,b]\otimes Hm,\\
& [a,a]\otimes m=a\otimes Pn
\end{align*}
where $a\otimes m$ is linear in $m$ and $[a,b]\otimes n$ is linear
in $a,b$ and $n$. For example, the quadratic $\mathbb Z$-modules
\begin{align*}
& \mathbb Z^\otimes=(\mathbb Z\buildrel{(1,1)}\over\rightarrow
\mathbb Z\oplus\mathbb Z
\buildrel{(1,1)}\over\rightarrow \mathbb Z)\\
& \mathbb Z^\Lambda=(0\rightarrow \mathbb Z\rightarrow 0)\\
& \mathbb Z^\Gamma=(\mathbb Z\buildrel{1}\over\rightarrow \mathbb
Z \buildrel{2}\over\rightarrow \mathbb Z)\\
& \mathbb Z^S=(\mathbb Z\buildrel{2}\over\rightarrow \mathbb Z
\buildrel{1}\over\rightarrow \mathbb Z)\\
& \mathbb Z_2=(\mathbb Z_2\rightarrow 0\rightarrow \mathbb Z_2)
\end{align*}
define the tensor square, exterior square, Whitehead's
$\Gamma$-functor, symmetric tensor square and the functor
$-\otimes \mathbb Z_2$ respectively, i.e. for every abelian group
$A$, one has
\begin{align*}
& A\otimes \mathbb Z^\otimes=A\otimes A,\ A\otimes \mathbb
Z^\Lambda=\Lambda^2A,\ A\otimes \mathbb Z^\Gamma=\Gamma A,\\
& A\otimes \mathbb Z^S=SP^2(A),\ A\otimes (\mathbb Z_2\to
0\to\mathbb Z_2)=A\otimes \mathbb Z_2
\end{align*}
where $SP^2(A)=A\otimes A/(a\otimes b-b\otimes a,\ a,b\in A)$.

Recall the definitions of some functors (see \cite{HTH}). Given
abelian groups $A,B$ we choose free resolutions
\begin{equation}\label{resolu}
d_A: X_1\to X_0,\ d_B: Y_1\to Y_0,
\end{equation}
of $A$ and $B$ respectively so that $d_A=C(A,0),\ d_B=C(B,0)$. For
a quadratic $\mathbb Z$-module $M$, one gets the chain complex
$M_\#(d_B)$ defined by
$$
Y_1\otimes Y_1\otimes M_{ee}\buildrel{\delta_2}\over\rightarrow
Y_1\otimes M\oplus Y_1\otimes Y_0\otimes
M_{ee}\buildrel{\delta_1}\over\rightarrow Y_0\otimes M
$$
with
\begin{align*}
& \delta_1(a\otimes m)=(d_Ba)\otimes m,\\
& \delta_1([a,a']\otimes n)=[d_Ba,a']\otimes n,\\
& \delta_1(a\otimes b\otimes n)=[d_Ba,b]\otimes n,\\
& \delta_2(a\otimes a'\otimes n)=-a\otimes d_Ba'\otimes
n+[a,d_Ba']\otimes n
\end{align*}
for $a,a'\in Y_1, b\in Y_0, m\in M_e, n\in M_{ee}$. Then
$coker(\delta_0)=A\otimes M$. Define the torsion functors by
\begin{align*}
& A*'M=ker(\partial_1)/im(\partial_2),\\
& A*''M=ker(\partial_2).
\end{align*}
Taking the torsion functors for $\mathbb Z^\Lambda$ and  $\mathbb
Z^\Gamma$ one, in fact, gets the Eilenberg-MacLane functors, see
6.2.9 and 6.2.10 in \cite{HTH}:
$$
A*'\mathbb Z^\Lambda=\Omega(A),\ A*'\mathbb Z^\Gamma=R(A).
$$

Define also the pseudo-analogs of the above torsion functors:
\begin{align*}
& \Lambda^2T_\#(A,B)=[d_A, \mathbb Z^\Lambda_\#d_B],\\
& \Gamma T_\#(A,B)=[d_A, \mathbb Z^\Gamma_\#d_B].
\end{align*}
Also recall the functor $L_\#(A,B)$, which for finitely generated
abelian groups $A$, $B$ is defined as
$$
L_\#(A,B)=B\otimes L(A).
$$
Here $L(A)$ is the quadratic $\mathbb Z$-module
$$
L(A)=(Hom(A,\mathbb Z_2)\buildrel{\partial}\over\rightarrow
Ext(A,\mathbb Z)\buildrel{0}\over\rightarrow Hom(A,\mathbb Z_2))
$$
where the map $\partial$ is the connecting homomorphism induced by
the exact sequence $\mathbb Z\buildrel{2}\over\rightarrow \mathbb
Z\to \mathbb Z_2$ (see 6.2.13 \cite{HTH} for the definition of
$L_\#(A,B)$ for arbitrary abelian groups). The bifunctors $\Lambda
T_\#, \Gamma T_\#$ and $L_\#$ are bype functors (see \cite{HTH})
in the sense that for abelian groups $A,B$ there are the following
binatural short exact sequences:
\begin{align*}
& 0\to Ext(A,\Omega(B))\to \Lambda T_\#(A,B)\to Hom(A,\Lambda^2(B))\to 0,\\
& 0\to Ext(A,R(B))\to \Gamma T_\#(A,B)\to Hom(A,\Gamma(B))\to 0,\\
& 0\to Ext(A,\Lambda^2(B))\to L_\#(A,B)\to Hom(A,B\otimes \mathbb
Z_2)\to 0.
\end{align*}
For every quadratic $\mathbb Z$-module and $n\geq 0$ there are
natural suspension maps
\begin{align*}
\Sigma^n: [d_A,M_\#d_B]\to [d_A[n],M_\#(d_B[n])]
\end{align*}
which in the case $M=\mathbb Z^\Lambda,\ M=\mathbb Z^\Gamma$
stabilize as
\begin{align*}
& \Sigma^n: \Lambda T_\#(A,B)\to Ext(A,B\otimes \mathbb Z_2),\ n>2\\
& \Sigma^n: \Gamma T_\#(A,B)\to Ext(A,B*\mathbb Z_2)\oplus
Hom(A,B\otimes \mathbb Z_2),\ n>3
\end{align*}

Let $d_C: Z_1\to Z_0$ be a short free resolution for an abelian
group $C$. Define the functor
$$
Trp(A;B,C)=[d_A, d_B\otimes d_C]=H_0(A,d_B\otimes d_C).
$$
Clearly one has the following short exact sequence which splits
(unnaturally)
$$
0\to Ext(A, B*C)\to Trp(A;B,C)\to Hom(A,B\otimes C)\to 1.
$$
The functor $Trp$ is an analogue of the triple torsion functor of
MacLane \cite{Maclane}. \vspace{.5cm}
\section{Quadratic bypes}\label{bypesec}
\vspace{.5cm} Let ${\sf Chain}_r$ be the category of $r$-reduced
chain complexes in ${\sf Ab}$. Let $M$ be a quadratic $\mathbb
Z$-module. Then we obtain a functor $(r\geq 0):$
$$
\begin{cases}
M_{\#}: {\sf Chain}_r\to {\sf Chain}_r\\
M_{\#}(Y)=N((N^{-1}Y)\otimes M),
\end{cases}
$$
where $\otimes M$ is the functor ${\sf Ab}\to {\sf Ab}$ given by
the quadratic tensor product.

An {\it $r$-reduced quadratic M-bype} is a triple $(Y,b,\beta)$
with the following properties. First $Y$ is a chain complex in
${\sf Chain}_r$. Let $B_n=H_nY$ be the homology of $Y$. Moreover
$b=\{b_n\}$ is a sequence of elements
$$
b_n\in Hom(B_n, H_{n-1}M_{\#}Y).
$$
Hence we get the homomorphism
$$
Ext(B_n,B_{n+1})\buildrel{(b_{n+1})_*}\over{\rightarrow} Ext(B_n,
H_nM_{\#}Y) \buildrel{\Delta}\over{\rightarrowtail}
H_{n-1}(B_n,M_{\#}Y),
$$
where the right hand side is the pseudo-homology. Finally
$\beta=\{\beta_n\}$ is a sequence of elements
$$
\beta_n\in H_{n-1}(B_n, M_{\#}(Y))/im(\Delta(b_{n+1})_*)
$$
with
$$
\mu\beta_n=b_n.
$$
A {\it morphism} between quadratic $M$-bypes
$$
\phi: (Y,b,\beta)\to (Y',b',\beta')
$$
is a morphism $\phi: Y\to Y'$ in ${\sf Ho(Chain)}_r$ with the
following property: the diagram
\begin{equation}\label{label1}
\begin{CD}
H_nY=B_n @>b_n>> H_{n-1}(M_{\#}Y)\\
@VV{\phi_*}V @VV{\phi_*}V\\
H_nY'=B_n' @>b_n'>> H_{n-1}(M_{\#}Y')
\end{CD}
\end{equation}
commutes for all $n$ and
\begin{equation}\label{label2}
\phi_*\beta_n=\phi^*\beta_n'
\end{equation}
in $H_{n-1}(B_n, M_{\#}Y')/\Delta(b_{n+1}')_*Ext(B_n,B_{n+1}').$
Let
$$ {\sf bype}_r(M)
$$
be the category of $r$-reduced quadratic $M$-bypes and such
morphisms.

Consider the quadratic $\mathbb Z$-module $\mathbb
Z^{\Lambda}=(0\to \mathbb Z\to 0)$ given by the exterior square
$\Lambda^2(A)=A\otimes \mathbb Z^{\Lambda}$.

\begin{theorem}\label{mainthe}
Let $r\geq 1$. Then there is a functor
$$
\lambda: {\sf Ho(sNil)}_r\to {\sf bype}_r(\mathbb Z^{\Lambda})
$$
which is representative and reflects isomorphisms. Moreover the
restriction of $\lambda$ to the subcategories of isomorphisms is a
full functor.
\end{theorem}
\begin{cor}
Let $r\geq 1$. Then the homotopy type of an $r$-reduced simplicial
2-nilpotent group $G$ is completely determined by the bype
$\lambda G$. In fact, the functor $\lambda$ induces a 1-1
correspondence between isomorphism types in ${\sf Ho(sNil)}_r$ and
isomorphism types in ${\sf bype}_r(\mathbb Z^\Lambda)$. Moreover
for each object $G$ in ${\sf (sNil)}_r$ the functor $\lambda$
induces a surjection of automorphism groups
$$
Aut_{{\sf Ho(sNil)}_r}(G)\buildrel{\lambda}\over\twoheadrightarrow
Aut_{{\sf bype}_r(\mathbb Z^\Lambda)}(\lambda G).
$$
\end{cor}
The functor $\lambda$ carries a free object $G$ in ${\sf sNil}$ to
$(Y,b,\beta)$ where $Y=NX$ with $X=ab(G)$. We may assume that $G$
is a free object. Then $G$ is part of the exact sequence
$$
\Lambda^2X\rightarrowtail G\to X
$$
in $\sf sNil$. This exact sequence is a fiber sequence in the
model category $\sf sNil$. Hence we get the connecting
homomorphisms $\partial=\partial_G$,
\begin{multline*}
b_n: B_n=H_nNX=\pi_nX\buildrel{\partial}\over{\rightarrow}
\pi_{n-1}\Lambda^2X=H_{n-1}N\Lambda^2X=H_{n-1}N\Lambda^2N^{-1}Y=H_{n-1}\mathbb
Z_{\#}^\Lambda Y,
\end{multline*}
where the isomorphism for $Y$ is given by $Y=NX$. This connecting
homomorphism defines the sequence of elements $b_n$ in
$\lambda(G)=(Y,b,\beta)$. We also can apply the functor
$\pi_n(A,-)$ of homotopy groups with coefficients in $A$ to the
fiber sequence above. This yields the following connecting
homomorphism $\partial=\partial_G$:

\begin{equation}
\xymatrix{Ext(A,\pi_{n+1}X)\ar[d]\ar[r]& \pi_n(A,X)\ar[ddd]^\partial\ar[r]^{\mu}&Hom(A,\pi_nX)\ar[d]^{\simeq}\\
Ext(A,B_{n+1})\ar[d]^{(b_{n+1})_*}&&Hom(A,B_n)\ar[d]^{(b_n)_*}\\
Ext(A,H_n\mathbb Z_{\#}^{\Lambda}Y)\ar[d]^\simeq &&Hom(A,H_{n-1}\mathbb Z_{\#}^{\Lambda}Y)\ar[d]^\simeq\\
Ext(A,\pi_n\Lambda^2X)\ar[r]^\Delta & \pi_{n-1}(A,\Lambda^2X)\ar[d]^\simeq\ar[r]^\mu &Hom(A,\pi_{n-1}\Lambda^2X)\\
&H_{n-1}(A,\mathbb Z_{\#}^\Lambda Y)&}
\end{equation}
In this diagram we set $A=B_n$ and define
$$
\beta_n=\{\partial\mu^{-1}(1_{B_n})\},
$$
where the right hand side is an element in the cokernel of
$\Delta(b_{n+1})_*$. It is clear that $\mu(\beta_n)=b_n$ and that
$\lambda$ is a well defined functor since the connecting
homomorphism $\partial$ is natural in $G$ and $A$. \vspace{.5cm}
\section{Homological quadratic bypes}
\vspace{.5cm} In this section $M$ is a quadratic $\mathbb
Z$-module. Then the category of $r$-reduced $M$-bypes admits
itself a detecting functor
$$
{\sf bype}_r(M)\buildrel{h}\over\rightarrow {\sf Hbype}_r(M),
$$
where the right hand side is the category of $r$-reduced
homological quadratic bypes which we define below. Let ${\sf
Chain}_r$ be the category of $r$-reduced chain complexes in ${\sf
Ab}$. We have the homology functor $H_*$
$$
H_*: {\sf Ho(Chain}_r) \to {\sf Ab}_r
$$
which admits a splitting functor $C$
$$
C:{\sf Ab}_r\to {\sf Ho(Chain}_r)
$$
defined by
$$
C(B)=\oplus_{n\in \mathbb Z}C(B_n,n),
$$
where $C(A,n)$ for an abelian group $A$ is the Moore chain
complex. Since $[C(A,n),C(A',n)]=Hom(A,A')$, we see that $C$ is a
well defined functor. Moreover for any $Y$ in ${\sf Chain}_r$ with
$B=H_*Y$ we can choose a weak equivalence
$$
C(B)\buildrel{\sim}\over\rightarrow Y.
$$
Recall that the {\it pseudo-homology} of a chain complex $Y$ in
${\sf Ab}_r$ is the set of homotopy classes of chain maps
$$
H_n(A,Y)=[C(A,n),Y].
$$
We have the short exact sequence
$$
0 \to Ext(A,H_{n+1}Y)\buildrel{\Delta}\over\rightarrow
H_n(A,Y)\buildrel{\mu}\over\rightarrow Hom(A,H_nY)\to 0.
$$
Using the quadratic $\mathbb Z$-module $M$ we define the composite
functor, see $\S 7$,
\begin{equation}\label{owo1}
Sq_n^M: {\sf Ab}_r\buildrel{C}\over\rightarrow{\sf Ho(Chain}_r)
\buildrel{M_\#}\over\rightarrow {\sf Ho(Chain}_r)
\buildrel{H_n}\over\rightarrow {\sf Ab}
\end{equation}
Moreover we define the bifunctor
\begin{equation}\label{owo2}
Sq_n^M: {\sf Ab}^{op}\times {\sf Ab}_r\to {\sf Ab}
\end{equation}
by the pseudo-homology $Sq_n^M(A,B)=H_n(A,M_{\#}CB).$ Hence we
have $Sq_n^M(B)=Sq_n^M(\mathbb Z,B)$ and one gets the binatural
short exact sequence
\begin{equation}\label{byper}
0\to Ext(A,Sq_{n+1}^M(B))\to Sq_n^M(A,B)\to Hom(A,Sq_n^M(B))\to 0.
\end{equation}
If $B\in {\sf Ab}_r$ is concentrated in one degree, that is,
$B=(D,m)$ with $B_m=D$ and $B_i=0$ for $i\neq m$, then we get the
functor
\begin{equation}
Sq_{n,m}^M: {\sf Ab}^{op}\times {\sf Ab}\to {\sf Ab}
\end{equation}
with $Sq_{n,m}^M(A,D)=Sq_n^M(A,(D,m))=H_n(A,M_\#C(D,m))$. Now we
set $Sq_{n,m}^M(D)=Sq_{n,m}^M(\mathbb Z,D)$ so we get the
binatural exact sequence
$$
0\to Ext(A,Sq_{n,m}^M(D))\to Sq_{n,m}^M(A,D)\to
Hom(A,Sq_{n,m}^M(D))\to 0
$$
as a special case of the exact sequence for (\ref{owo2}) above.
Here we have as an example the homotopy groups of a Moore objects
$M(D,n)_{nil}$ in ${\sf sNil}$ given by $(k\geq 1)$
$$
\pi_{n+k}M(D,n)_{nil}=Sq_{n+k,n}^{\mathbb Z^\Lambda}(D)
$$
which was computed in theorem \ref{mooreh}.

Since $M_{\#}$ is a quadratic functor we get morphisms in ${\sf
Chain}_r$
$$
M_{\#}(Y)\buildrel{H}\over\rightarrow
M_\#(Y|Y)\buildrel{P}\over\rightarrow M_\#(Y)
$$
If $Y$ is free as an $R$-module we have by the Eilenberg-Zilber
theorem
$$
M_\#(Y|Y)=N((N^{-1}Y)\otimes (N^{-1}Y)\otimes M_{ee})\simeq
Y\otimes Y\otimes M_{ee}
$$
Hence we also get
$$
M_\#(Y)\buildrel{H}\over\rightarrow Y\otimes Y\otimes
M_{ee}\buildrel{P}\over\rightarrow M_\#(Y)
$$
by $H$ and $P$ above. For $M=\mathbb Z^\Lambda$ we compute the
functors $Sq_n^M$ explicitly below. We use the functors $Sq_n^M$
for the definition of the following category ${\sf Hbype}_r(M)$.

Let $r\geq 0$. An {\it $r$-reduced homological quatratic M-bype}
is a triple $(B,b,\beta)$ with the following properties. First $B$
is a graded abelian group in ${\sf Ab}_r$ and $b=\{b_n\}$ is a
sequence of elements
$$
b_n\in Hom(B_n, Sq_{n-1}^M(B)).
$$
Hence we get the homomorphism
$$
Ext(B_n,B_{n+1})\buildrel{(b_{n+1})_*}\over\rightarrow Ext(B_n,
Sq_n^M(B)) \buildrel{\Delta}\over\rightarrowtail
Sq_{n-1}^M(B_n,B),
$$
where the right hand side is the pseudo-homology above. Now
$\beta=\{\beta_n\}$ is a sequence of elements
$$
\beta_n\in Sq_{n-1}^M(B_n,B)/im(\Delta(b_{n+1})_*)
$$
with
$$
\mu\beta_n=b_n.
$$
A {\it morphism} between quadratic $M$-bypes
$$
\phi: (B,b,\beta)\to (B',b',\beta')
$$
is a morphism $\phi: B\to B'$ in ${\sf Ab}_r$ with the following
property. There exists a chain map $\alpha: C(B)\to C(B')$ with
$H_*(\alpha)=0$ such that the chain map
$$
\phi_*+\alpha: M_\#CB\to M_\#CB',
$$
given by
$$
\phi_*+\alpha:=M_\#C(\phi)+P(\alpha\otimes C(\phi))H+M_\#(\alpha)
$$
satisfies (\ref{refg1}) and (\ref{refg2}): the diagram
\begin{equation}\label{refg1}
\begin{CD}
B_n @>{b_n}>> H_{n-1}(M_\#CB)=Sq_{n-1}(B)\\
@V{\phi}VV @V{(\phi_*+\alpha)_*}VV\\
B_n' @>{b_n'}>> H_{n-1}(M_\#CB')=Sq_{n-1}(B')
\end{CD}
\end{equation}
commutes for all $n$ and
\begin{equation}\label{refg2}
(\phi_*+\alpha)_*\beta_n=\phi^*\beta_n'
\end{equation}
in $Sq_{n-1}(B_n,B')/\Delta(b_{n+1}')_*Ext(B_n,B_{n+1}').$ Diagram
(\ref{refg1}) shows that $(\phi_*+\alpha)_*$ in (\ref{refg2}) is
well defined by $H_{n-1}(B_n,\phi_*+\alpha)$. Let
$$
{\sf Hbype}_r(M)
$$
be the category of $r$-reduced homological quadratic $M$-bypes and
such morphisms. We point out that the homology class of a chain
map $\alpha: CB\to CB'$ with $H_*\alpha=0$ is given by a sequence
of elements $$ \alpha_n\in Ext(B_n,B_{n+1}'),\ n\in\mathbb Z,
$$
and the induced maps
\begin{align*}
& \alpha_*: H_{n-1}(M_\#CB)\to H_{n-1}(M_\#CB')\\
& \alpha_*: H_{n-1}(B_n, M_\#CB)\to H_{n-1}(B_n, M_\#CB')
\end{align*}
are given by the maps induced by $Ext$ in section \ref{section7}.
These induced maps are needed for the computation of
$(\phi_*+\alpha)_*$ in (\ref{refg1}) and (\ref{refg2}) above.
\begin{prop}\label{propp2}
There is a functor
$$
h: {\sf bype}_r(M)\to {\sf Hbype}_r(M),
$$
which carries $(Y,b,\beta)$ to $(H_*Y,b,\beta)$. The functor $h$
reflects isomorphisms and is representative and full.
\end{prop}
Here we use the equivalence $CB\buildrel{\sim}\over\rightarrow Y$
with $B=H_*Y$. Combining Proposition \ref{propp2} and Theorem
\ref{mainthe} we get the following main result of this paper.

\begin{theorem}
Let $r\geq 1$. Then there is a functor
$$
h\lambda: {\sf Ho(sNil)}_r\to {\sf Hbype}_r(\mathbb Z^\Lambda)
$$
which reflects isomorphisms and is representative. Moreover the
restriction of $h\lambda$ to the subcategories of isomorphisms is
a full functor.
\end{theorem}

In the next section we compute $Sq_n^{\mathbb Z^\Lambda}$ needed
in the definition of ${\sf Hbype}_r(\mathbb Z^\Lambda)$.
\vspace{.5cm}
\section{Computation of $Sq_n^M$ for $M=\mathbb Z^\Lambda$}
\vspace{.5cm} For an abelian group $A$ and for a graded abelian
group $B$ in ${\sf Ab}_r$ and $M=\mathbb Z^\Lambda$ we compute
$Sq_n^M(A,B)$ with $Sq_n^M(B)=Sq_n^M(\mathbb Z,B)$ as follows.

Given abelian groups $A,D,E$, denote
\begin{align*} & Sq_{n,m}^{\mathbb Z^\Lambda}(A,D)=[C(A,n),\mathbb Z^{\Lambda}_\#C(D,m)],\ n,m\geq
1,\ \text{see}\ 8.4\\
& Sq_{n,i,j}^{\mathbb Z^\Lambda}(A;D,E)=[C(A,n),C(D,i)\otimes
C(E,j)],\ n\geq 1,\ 1\leq i<j
\end{align*} Then for $B\in {\sf Ab}_r$, we have the following direct sum decomposition:
\begin{align*}
Sq_n^{\mathbb Z^\Lambda}(A,B)& =[C(A,n), \bigoplus_m\mathbb
Z^\Lambda_\#C(B_m,m)\oplus\bigoplus_{i<j}C(B_i,i)\otimes
C(B_j,j)]\\ & =\bigoplus_mSq_{n,m}^{\mathbb
Z^\Lambda}(A,B_m)\oplus \bigoplus_{i<j}Sq_{n,i,j}^{\mathbb
Z^\Lambda}(A;B_i,B_j).
\end{align*}
It is easy to see that for $i<j$,
$$
Sq_{n,i,j}^{\mathbb Z^\Lambda}(A;B_i,B_j)=\begin{cases}
Hom(A,B_i*B_j),\
n=i+j+1,\\
Ext(A,B_i\otimes B_j),\ n=i+j-1,\\
Trp(A;B_i,B_j),\ n=i+j,\\
0,\ \text{otherwise}
\end{cases}
$$

We describe the functors $Sq_{n,m}^{\mathbb Z^\Lambda}$ in the
following theorem.
\begin{theorem}\label{lamq}
\begin{align*}
Sq_{m+k,m}^{\mathbb Z^\Lambda}(A,D)=\begin{cases}
Ext(A,\Gamma(D)),\ k=0,\ m=1 & (I)\\
Ext(A,B\otimes \mathbb Z_2),\ k=0,\ m>1 & (II)\\
Hom(A,R(D)),\ k\ \text{even},\ k=m+1 & (III)\\
Hom(A,\Omega(D)),\ k\ \text{odd},\ k=m+1 & (IV)\\
Ext(A,D\otimes \mathbb Z_2)\oplus Hom(A,D*\mathbb
Z_2),\ k\ \text{even},\ 0<k<m-1 & (V)\\
Ext(A,D*\mathbb Z_2)\oplus Hom(A,D\otimes \mathbb Z_2),\ k\
\text{odd},\ 0<k<m-1 & (VI)\\
Ext(A,\Gamma(D))\oplus Hom(A,D*\mathbb Z_2),\
k\ \text{even},\ k=m-1 & (VII)\\
L_\#(A,D)\oplus Ext(A,D*\mathbb Z_2),\ k\ \text{odd},\
k=m-1 & (VIII)\\
\Lambda^2T_\#(A,D)\oplus Hom(A,D*\mathbb Z_2), k\ \text{even},\ k=m & (IX)\\
\Gamma T_\#(A,D),\ k\ \text{odd},\ k=m & (X)\\
0,\ \text{otherwise} & (XI)
\end{cases}
\end{align*}
\end{theorem}
We define the {\it stable} operator $Sq_k^{stable}$ by
$$
Sq_k^{stable}(A,D)=\begin{cases} Ext(A,D\otimes \mathbb Z_2)\oplus
Hom(A,D*\mathbb Z_2),\ k\ \text{even}>0,\\
Ext(A,D*\mathbb Z_2)\oplus Hom(A,D\otimes\mathbb Z_2),\ k\
\text{odd}>0
\end{cases}
$$
Then there is a canonical {\it stabilization map}
$$
\Sigma^\infty: Sq_{m+k,m}^{\mathbb Z^\Lambda}(A,D)\to
Sq_k^{stable}(A,D)
$$
which is binatural in $A,D$ and which is the identity for
$k<m-1$.\\

We collect the computations of $Sq_{n,m}^{\mathbb Z^\Lambda}(A,D)$
for low dimensions in the following tables:
\begin{align*}
& \begin{tabular}{ccccccc}
 $m\setminus n$ & \vline & 1 & 2 & 3 \\ \hline
1 & \vline & $Ext(A,\Gamma(D))$ & $\Gamma T_\#(A,D)$ & $Hom(A,R(D))$\\
2 & \vline & 0 & $Ext(A,D\otimes \mathbb Z_2)$ & $L_\#(A,D)\oplus
Ext(A,D*\mathbb Z_2)$\\
3 & \vline & 0 & 0 & $Ext(A,D\otimes \mathbb Z_2)$\\
4 & \vline & 0 & 0 & 0\\
\end{tabular}\\
&
\begin{tabular}{ccccccc}
$m\setminus n$ & \vline & 4 & 5\\
\hline 2 & \vline & $\Lambda^2T_\#(A,D)\oplus Hom(A,D*\mathbb
Z_2)$ &
$Hom(A,\Omega(D))$ \\
3 & \vline & $Ext(A,D*\mathbb Z_2)\oplus Hom(A,D\oplus \mathbb
Z_2)$ &
$Ext(A,\Gamma(D))\oplus Hom(A,D*\mathbb Z_2)$ \\
4 & \vline & $Ext(A,D\otimes \mathbb Z_2)$ & $Ext(A,D\otimes
\mathbb Z_2)\oplus Hom(A,D*\mathbb Z_2)$ \\
5 & \vline & 0 & $Ext(A,D\otimes\mathbb Z_2)$\\
\end{tabular}\\
&
\begin{tabular}{ccccccc}
 $m\setminus n$ & \vline & 6 & 7\\ \hline
2 & \vline & 0 & 0\\
3 & \vline & $\Gamma T_\#(A,D)$ & $Hom(A,R(D))$\\
4 & \vline & $Ext(A,D\otimes\mathbb Z_2)\oplus
Hom(A,D*\mathbb Z_2)$ & $L_\#(A,D)\oplus Ext(A,D*\mathbb Z_2)$\\
5 & \vline & $Ext(A,D*\mathbb Z_2)\oplus Hom(A,D\otimes \mathbb
Z_2)$ &
$Ext(A,D*\mathbb Z_2)\oplus Hom(A,D\otimes \mathbb Z_2)$\\
6 & \vline & $Ext(A,D\otimes\mathbb Z_2)$ & $Ext(A,D\otimes\mathbb
Z_2)\oplus
Hom(A,D*\mathbb Z_2)$\\
\end{tabular}
\end{align*}

\vspace{.5cm} \noindent{\bf Example.} If $B\in {\sf Ab}_1$ is
concentrated in degree $1,2,3$ with $B_3$ free abelian then the
only invariants of $(B,b,\beta)$ in ${\sf Hbype}_1(\mathbb
Z^\Lambda)$ are given by
$$
b_3:H_3\to \Gamma(H_2)\ \text{and}\ \beta\in Ext(H_3,coker(b_3))
$$
and they classify 1-connected 4-dimensional homotopy types by a
classical result of J.H.C. Whitehead \cite{Whitehead}, see 3.5.6
\cite{HTH}.\\

\noindent{\bf Example.} Let $G$ be a reduced 2-nilpotent
simplicial group so that $NG$ is a rational vector space. Then for
$B=H_*NG$ the only invariants in $(B,b,\beta)=h\lambda(G)$ are the
homomorphisms
$$
b_n:B_n\to Sq_{n-1}^{\mathbb Z^\Lambda}(B)=[B,B]_{n-1}
$$
where the Lie bracket $[\ ,\ ]$ in the free Lie algebra $L(B)$
satisfies $[x,y]=-(-1)^{|x||y|}[y,x]$ so that
$$
[B_i,B_i]=\begin{cases} SP^2(B_i)=\Gamma(B_i)\ \text{for}\ n\
\text{odd}\\
\Lambda^2(B_i)\ \text{for}\ n\ \text{even}
\end{cases}
$$
Here we use the isomorphism $\mathbb Z^\Gamma\otimes \mathbb
Q=\mathbb Z^S\otimes \mathbb Q$, compare section 6. The invariant
$b_n$ coincides with the differential in the 2-nilpotent
differential Lie algebra associated to $G$ in the work of Quillen
\cite{QuillenRational}.

\vspace{.5cm}
\section{Maps induced by $Ext(A,B)$}\label{section7}
\vspace{.5cm} Given abelian groups $A,B$ and their resolutions
(\ref{resolu}), consider an element $\alpha\in Ext(A,B)$. This
element can be represented as a certain diagram of the form
$\alpha_*: d_A\to d_B[1]$:
$$
\xyma{& X_1
\ar@{->}[r]^{d_A} \ar@{->}[d]^{\alpha_*} & X_0\\
Y_1 \ar@{->}[r]^{d_B} & Y_0 & &\\
}$$ For a given $n\geq 1,$ consider the shifting
$\alpha_*[n]:C(A,n)\to C(B,n+1)$ with $C(A,n)=d_A[n]$, which
defines the map of simplicial abelian groups
$$N^{-1}\alpha_*[n]:N^{-1}C(A,n)\to N^{-1}C(B,n+1)
$$
and hence, for a given quadratic $\mathbb Z$-module $M$, there is
a map
\begin{equation}\label{rd}
\xyma{N((N^{-1}C(A,n))\otimes M)
\ar@{->}[r]^{M_\#\alpha_*[n]} \ar@{=}[d] & N((N^{-1}C(B,n+1))\otimes M) \ar@{=}[d]\\
M_\#C(A,n) \ar@{->}[r] & M_\#C(B,n+1)\\
}\end{equation} Then, for every $k\geq 0$ this diagram induces for
abelian groups $D,A,B$ the natural quadratic maps
\begin{align*}&\ ^k[n]_M: Ext(A,B)\to Hom(Sq_{n+k,n}^M(A), Sq_{n+k,n+1}^M(B))\ \text{and}\\
& ^k[n]_M: Ext(A,B)\to Hom(Sq_{n+1}^M(D,A), Sq_{n+k,n+1}^M(D,B)).
\end{align*}
The description of the homotopy groups of Moore space in {\sf
sNil} (see theorem \ref{mooreh} or theorem 9.1) implies that for
the case $M=\mathbb Z^\Lambda$, the complete list of maps $^k[n]$
is the following: \begin{align*} & ^2[1]_{\mathbb Z^\Lambda}:
Ext(A,B)\to
Hom(R(A),B\otimes \mathbb Z_2),\\
& ^2[2]_{\mathbb Z^\Lambda}: Ext(A,B)\to Hom(\Lambda^2(A)\oplus
A*\mathbb Z_2,
B\otimes \mathbb Z_2),\\
& ^2[3]_{\mathbb Z^\Lambda}: Ext(A,B)\to Hom(A*\mathbb Z_2,
B\otimes \mathbb
Z_2),\\
& ^3[2]_{\mathbb Z^\Lambda}: Ext(A,B)\to Hom(\Omega(A),B*\mathbb Z_2),\\
& ^3[3]_{\mathbb Z^\Lambda}: Ext(A,B)\to Hom(\Gamma(A), B*\mathbb Z_2),\\
& ^3[4]_{\mathbb Z^\Lambda}: Ext(A,B)\to Hom(A\otimes \mathbb Z_2,
B*\mathbb Z_2)
\end{align*}
Let
\begin{equation}\label{mna}0\to B\to E\to A\to 0\end{equation} be
a short exact sequence, which presents an element $\alpha\in
Ext(A,B)$. Applying the functor $-\otimes \mathbb Z_2$, we get the
following long exact sequence
$$
0\to B*\mathbb Z_2\to E*\mathbb Z_2\to A*\mathbb
Z_2\buildrel{\partial(\alpha)}\over\rightarrow B\otimes\mathbb
Z_2\to E\otimes\mathbb Z_2\to A\otimes\mathbb Z_2\to 0
$$
and the map $^2[3]_{\mathbb Z^\Lambda}$ is given by setting
$\alpha\mapsto
\partial(\alpha)\in Hom(A*\mathbb Z_2,B\otimes \mathbb Z_2)$.
\begin{prop}\label{pravesh}
The maps $^3[2]_{\mathbb Z^\Lambda}$, $^3[3]_{\mathbb Z^\Lambda}$
and $^3[4]_{\mathbb Z^\Lambda}$ are zero maps, $^2[3]_{\mathbb
Z^\Lambda}$ is given by setting $\alpha\mapsto \partial(\alpha),$
the maps $^2[1]_{\mathbb Z^\Lambda}$ and $^2[2]_{\mathbb
Z^\Lambda}$ are induced by the natural maps $R(A)\to A*\mathbb
Z_2,$ $\Lambda^2(A)\oplus A*\mathbb Z_2\to A*\mathbb Z_2$ and the
map $\partial(\alpha)$.
\end{prop}
\begin{proof}
The exact sequence of quadratic $\mathbb Z$-modules
$$
0\to \mathbb Z^\Gamma\to \mathbb Z^\otimes\to \mathbb Z^\Lambda\to
0
$$
and the epimorphism of quadratic $\mathbb Z$-modules
$$
\mathbb Z^\Gamma\to \mathbb Z_2
$$
induce the following natural commutative diagram
\begin{equation}\label{drf}
\begin{CD}
H_{n+k}\mathbb Z^\Lambda_\#C(A,n) @>>> H_{n+k-1}\mathbb
Z^\Gamma_\#C(A,n) @>>> H_{n+k-1}C(A,n)\otimes \mathbb Z_2\\
@V{^k[n]_{\mathbb Z^\Lambda}}VV @V{^{k-1}[n]_{\mathbb Z^\Gamma}}VV @V{^{k-1}[n]_{\mathbb Z_2}}VV\\
H_{n+k}\mathbb Z^\Lambda_\#C(B,n+1) @>>> H_{n+k-1}\mathbb
Z^\Gamma_\#C(B,n+1) @>>> H_{n+k-1} C(B,n+1)\otimes\mathbb Z_2
\end{CD}
\end{equation}
For $n>2$ and $k=1,2$, the diagram (\ref{drf}) has the following
structure:
\begin{align*}
& k=1 \xyma{& A\otimes \mathbb Z_2 \ar@{=}[r]
\ar@{->}[d]^{^1[n]_{\mathbb Z^\Lambda}} & A\otimes\mathbb Z_2
\ar@{=}[r]
\ar@{->}[d]^{^0[n]_{\mathbb Z^\Gamma}} & A\otimes\mathbb Z_2 \ar@{->}[d]^{^0[n]_{\mathbb Z_2}}\\
& B*\mathbb Z_2 \ar@{=}[r] & B*\mathbb Z_2 \ar@{->}[r] & 0 \\
}\\ & k=2 \xyma{& A* \mathbb Z_2 \ar@{=}[r]
\ar@{->}[d]^{^2[n]_{\mathbb Z^\Lambda}} & A*\mathbb Z_2 \ar@{=}[r]
\ar@{->}[d]^{^1[n]_{\mathbb Z^\Gamma}}
& A*\mathbb Z_2 \ar@{->}[d]^{^1[n]_{\mathbb Z_2}}\\
& B\otimes \mathbb Z_2 \ar@{=}[r] & B\otimes \mathbb Z_2 \ar@{=}[r] & B\otimes \mathbb Z_2 \\
}
\end{align*}
where, clearly, $^1[n]_{\mathbb Z_2}=\partial(\alpha)$. Now,
taking the suspension functors, we get the following diagram with
commutative squares
$$
\begin{CD}
H_3(\mathbb Z_\#^\Lambda C(A,1)) @>{\Sigma}>> H_4(\mathbb
Z_\#^\Lambda C(A,2)) @>{\Sigma}>> H_5(\mathbb Z_\#^\Lambda
C(A,3)\\
@V{ ^2[1]_{\mathbb Z^\Lambda}}VV @V{ ^2[2]_{\mathbb Z^\Lambda}}VV @V{^2[3]_{\mathbb Z^\Lambda}}VV\\
H_3(\mathbb Z_\#^\Lambda C(B,2)) @>{\Sigma}>> H_4(\mathbb
Z_\#^\Lambda C(B,3)) @>{\Sigma}>> H_5(\mathbb Z_\#^\Lambda C(B,4)
\end{CD}
$$
which is
$$
\begin{CD}
R(A) @>{\Sigma}>> \Lambda^2(A)\oplus A*\mathbb Z_2 @>{\Sigma}>>
A*\mathbb Z_2\\
@V{ ^2[1]_{\mathbb Z^\Lambda}}VV @V{ ^2[2]_{\mathbb Z^\Lambda}}VV @V{\partial(\alpha)}VV\\
B\otimes \mathbb Z_2 @= B\otimes \mathbb Z_2 @= B\otimes \mathbb
Z_2
\end{CD}
$$
Analogically we get the suspension diagram
$$
\begin{CD}
\Omega(A) @>{\Sigma}>> \Gamma(A) @>{\Sigma}>>
A\otimes Z_2\\
@V{ ^3[2]_{\mathbb Z^\Lambda}}VV @V{ ^3[3]_{\mathbb Z^\Lambda}}VV
@V{ ^3[4]_{\mathbb Z^\Lambda}}V{0}V\\
B* \mathbb Z_2 @= B* \mathbb Z_2 @= B* \mathbb Z_2
\end{CD}
$$
\end{proof}

\noindent{\bf Remark.} Proposition \ref{pravesh} implies that the
maps $^k[n]_{\mathbb Z^\Lambda}$ are zero for all $n$ and odd
$k$.\\

Given abelian groups $A,D,D'$ and an element $\alpha\in
Ext(D,D')=[C(D,0),C(D',1)]$ the shifting map $\alpha_*[m]:
C(D,m)\to C(D',m+1)$ induces the map of pseudo-homology
$$
\{n,m\}_M: Sq_{n,m}^M(A,D)\to Sq_{n,m+1}^M(A,D')
$$
The description of these maps follows from Proposition
\ref{pravesh}. We have the following commutative diagram
\begin{equation}\label{nikas}
\xyma{Ext(A,Sq_{n+1,m}^{\mathbb Z^\Lambda}(D)) \ar@{>->}[r]
\ar@{->}[d]^{Ext(A,\ ^{n+1-m}[m]_{\mathbb Z^\Lambda})} &
Sq_{n,m}^{\mathbb Z^\Lambda}(A,D) \ar@{->>}[r]
\ar@{->}[d]^{\{n,m\}_{\mathbb Z^\Lambda}} &
Hom(A,Sq_{n,m}^{\mathbb Z^\Lambda}(D)) \ar@{->}[d]^{Hom(A,\ ^{n-m}[m]_{\mathbb Z^\Lambda})}\\
Ext(A,Sq_{n+1,m+1}^{\mathbb Z^\Lambda}(D')) \ar@{>->}[r] &
Sq_{n,m+1}^{\mathbb Z^\Lambda}(A,D') \ar@{->>}[r] & Hom(A,
Sq_{n,m+1}^{\mathbb Z^\Lambda}(D'))
\\
}
\end{equation}
Since $^k[m]_{\mathbb Z^\Lambda}$ is the zero map for odd $k$, one
of the maps either $Ext(A,\ ^{n+1-m}[m]_{\mathbb Z^\Lambda})$ or
$Hom(A,\ ^{n-m}[m]_{\mathbb Z^\Lambda})$ is zero and the map
$\{n,m\}_{\mathbb Z^\Lambda}$ is defined via diagram
(\ref{nikas}). \vspace{.5cm}
\section{Homotopy types of spectra in
{\sf sNil} and $\mathcal F$-modules.}\vspace{.5cm} The homotopy
theory $\sf specnil$ of connective spectra in the model category
${\sf sNil}$ is defined in \cite{SMC}. For homotopy categories we
have the stabilization functor
$$
\Sigma^\infty: {\sf Ho(sNil)\to Ho(specnil)}
$$
which carries a free object $G$ in ${\sf sNil}$ to its suspension
spectrum. This functor has an analogue on the level of homological
bypes in the sense that one has the diagram of functors
\begin{equation}\label{lifafa}
\begin{CD} {\sf Ho(sNil)}_r @>{\Sigma^\infty}>> {\sf
Ho(specnil)}\\
@V{\lambda}VV @VV{\lambda^{spec}}V\\
{\sf Hbype}_r(\mathbb Z^\Lambda) @>{\Sigma^\infty}>> {\sf
mod}(\mathcal F,\Delta,\mu)
\end{CD}
\end{equation}
which commutes up to a natural isomorphism. Here the category
${\sf mod}(\mathcal F,\Delta,\mu)$ is defined as follows. Let
$$
\mathcal F=T_{\mathbb Z_2}(Sq_2^{nil},Sq_4^{nil},\dots)
$$
be the free graded $\mathbb Z_2$-algebra generated by elements
$Sq_k^{nil}$ of degree $-k$ for $k$ even $>0$. An object
$(H,H(2))$ in ${\sf mod}(\mathcal F,\Delta, \mu)$ is given by a
graded abelian group $H$ with $H_n=0,$ for $n<0$ and a graded
$\mathbb Z_2$-vector space $H(2)$ together with a short exact
sequence
$$
0\to H\otimes \mathbb Z_2\buildrel{\Delta}\over\rightarrow
H(2)\buildrel{\mu}\over\rightarrow (H*\mathbb Z_2)[-1]\to 0
$$
where in addition $H(2)$ is an $\mathcal F$-module, i.e. maps
$$
Sq_k^{nil}: H(2)_n\to H(2)_{n-k}
$$
are defined for $k$ even $>0$. A morphism $\phi: (H,H(2))\to (\bar
H,\bar H(2))$ is a homomorphism $\phi: H\to \bar H$ between graded
abelian groups for which there exists a commutative diagram
$$
\begin{CD}
0 @>>> H\otimes\mathbb Z_2 @>{\Delta}>> H(2) @>{\mu}>> (H*\mathbb
Z_2)[-1] @>>> 0\\
@. @V{\phi\otimes 1}VV @V{\psi}VV @V{\phi*\mathbb Z_2}VV @.\\
0 @>>> \bar H\otimes \mathbb Z_2 @>{\Delta}>> \bar H(2) @>{\mu}>>
(\bar H*\mathbb Z_2)[-1] @>>> 0
\end{CD}
$$
where $\psi$ is a map of $\mathcal F$-modules.
\begin{theorem}\label{chaya}
There exists a functor $\lambda^{spec}$ for which diagram
(\ref{lifafa}) commutes up to natural isomorphism. Moreover
$\lambda^{spec}$ is representative and reflects isomorphisms and
the restriction of $\lambda^{spec}$ to the subcategories of
isomorphisms is a full functor.
\end{theorem}
The theorem shows that homotopy types of connected spectra in
${\sf sNil}$ are completely determined by an isomorphism type of
an $\mathcal F$-module $(H,H(2))$ in the category ${\sf
mod}(\mathcal F,\Delta, \mu)$. The $\mathcal F$-module structure
is related  to the action of the Steenrod algebra as follows. \\

\noindent{\bf Remark.} Let $G$ be an object in $({\sf sNil})_r$
which is given by the nilization of a space $X$. Then the
$\mathcal F$-module
$$
\Sigma^\infty \lambda(G)\cong \lambda^{spec}\Sigma^\infty
(G)=(H,H(2))
$$
is defined such that the operator $Sq_k^{nil}$ on $H(2)$ fits into
the following commutative diagram
$$
\begin{CD}
H(2)_n @>{Sq_k^{nil}}>> H(2)_{n-k}\\
@| @|\\
H_{n+1}(X,\mathbb Z_2) @>{(\chi Sq^k)}_*>> H_{n+1-k}(X,\mathbb
Z_2)
\end{CD}
$$
Here $Sq^k$ denotes the Steenrod square and $\chi$ is the
anti-isomorphism of the Steenrod algebra and $(\chi Sq^k)_*$ is
obtained by dualization. The commutativity of the diagram follows
from 8.13 \cite{B}.\\

\noindent{\it Proof of Theorem \ref{chaya}} We define the functor
$\Sigma^\infty$ in the bottom row of (\ref{lifafa}) by the
composite of functors
\begin{equation}\label{loha}
\Sigma^\infty: {\sf Hbype}_r(\mathbb Z^\Lambda)
\buildrel{\Sigma}\over\to {\sf Hbype}^\infty
\buildrel{\Theta}\over\cong {\sf mod}(\mathcal F,\Delta,\mu)
\end{equation}
where $\Theta$ is an isomorphism of categories.

For this we introduce the category ${\sf Hbype}^\infty$ of stable
homological bypes as follows. A {\it stable homological bype}
$(B,b,\beta)$ consists of a graded abelian group $B$ with $B_n=0$
for $n<0$ and homomorphisms $(n,k\in \mathbb Z,\ k\geq 1)$
\begin{align*}
& b_n^k: B_n\otimes\mathbb Z_2\to \begin{cases} B_{n-k}\otimes
\mathbb Z_2\ \text{for}\ $k$\ \text{even}\\ B_{n-k}*\mathbb Z_2\
\text{for}\ k\ \text{odd}
\end{cases}\\
& \beta_n^k: B_n*\mathbb Z_2\to \begin{cases} B_{n-k}* \mathbb
Z_2\ \text{for}\ $k$\ \text{even}\\ B_{n-k}\otimes\mathbb Z_2\
\text{for}\ k\ \text{odd}
\end{cases}
\end{align*}
Here $(B,b,\beta)$ is {\it equivalent} to $(\tilde B, \tilde
b,\tilde \beta)$ if $B=\tilde B$, $b=\tilde b$ and $\beta\sim
\tilde \beta$ in the sense that
$$
\tilde \beta_{n-1}^k=\beta_{n-1}^k+b_n^{k+1}\delta_n
$$
for some homomorphism $\delta_n: B_{n-1}*\mathbb Z_2\to B_n\otimes
\mathbb Z_2$. Objects of the category ${\sf Hbype}^\infty$ are
such equivalence classes $\{B,b,\beta\}$. A morphism
$\{B,b,\beta\}\to \{\bar B,\bar b,\bar\beta\}$ in ${\sf
Hbype}^\infty$ is a homomorpism $\phi: B\to \bar B$ of graded
abelian groups for which there exist $(j\in \mathbb Z)$
$$
\alpha_j\in Ext(B_j,\bar B_{j+1}\otimes \mathbb
Z_2)=Hom(B_j*\mathbb Z_2,\bar B_{j+1}\otimes \mathbb Z_2)
$$
such that $(k\ \text{even},\ n\in \mathbb Z)$
\begin{align*}
& \phi_*b_n^k+\alpha_{n-k}b_n^{k+1}=\bar b_n^k\phi_*\\
& \phi_*\beta_{n-1}^{k-1}+\alpha_{n-k}\beta_{n-1}^k\sim
\bar\beta_{n-1}^{k-1}\phi_*
\end{align*}
One can check that stabilization of $Sq_{m+k,m}(A,B)$ yields a
canonical functor $\Sigma$ in (\ref{loha}) which carries
$(B,b,\beta)$ to $(B,\Sigma^\infty b,\Sigma^\infty \beta)$.
Moreover the isomorphism $\Theta$ is given by setting
$$
\Theta(B,b,\beta)=(H,H(2))
$$
where $H=B$ and $H(2)=B_n\otimes \mathbb Z_2\oplus B_{n-1}*\mathbb
Z_2$ with $Sq_k^{nil}$ given by the matrix
$$
Sq_k^{nil}=\begin{pmatrix} b_n^k & \beta_{n-1}^{k-1}\\
b_n^{k+1} & \beta_{n-1}^k
\end{pmatrix}
$$
Here we choose $\beta$ in the equivalence class $\{B,b,\beta\}.$
One can check that $\Theta$ is a well defined isomorphism of
categories. In the stable range of ${\sf bype}_r(\mathbb
Z^\Lambda)$ the functor $\Sigma^\infty$ in (\ref{loha}) is a full
embedding of categories for all $r\geq 1$. This implies the result
on ${\sf specnil}$ in theorem \ref{chaya}.\ $\Box$
 \vspace{.5cm}
\section{Proof of Theorem \ref{mainthe}}
\vspace{.5cm} The proof is based on the theory of boundary
invariants developed in \cite{HTH} for CW-complexes. A similar
theory is available for CW-objects in the category ${\sf Nil}$
which are the free objects in section 4.

It is a well-known result of Kan \cite{HT} that the homotopy
theory in $({\sf sGr})_r$ is equivalent to the homotopy theory of
CW-complexes $X$ with trivial $r$-skeleton. In particular, the
generators of a free simplicial group $H$ correspond to the cells
of a CW-complex $X$ with $X\simeq B|H|$. This way we can associate
the boundary invariants in \cite{HTH} for the CW-complex $X$ to
the simplicial group $H$. If $G=nil(H)$ is the nilization of $G$
we get the commutative diagram of short exact sequences
$$
\begin{CD}
[H,H] @>>> H @>>> ab H\\
@VVV @VVV @|\\
\Lambda^2(abG) @>>> G @>>> ab G
\end{CD}
$$
where $[H,H]$ is the commutator subgroup. Hence we have natural
maps
\begin{align}
& \Gamma_{n+1}X=\pi_n[H,H]\to\pi_n(\Lambda^2(abG))\label{u1}\\
& \Gamma_n(A,X)=\pi_{n-1}(A,
[H,H])\to\pi_{n-1}(A,\Lambda^2(abG))\label{u2}.
\end{align}
Here the groups $\Gamma_{n+1}X$ and $\Gamma_n(A,X)$ are defined in
\cite{HTH}, compare also section 2.4 in \cite{HTH}. The boundary
invariants of a space $X$ are developed in \cite{HTH}. In a
similar way one gets the boundary invariants of a free object in
${\sf sNil}$ such that the natural maps (\ref{u1}) and (\ref{u2})
carry boundary invariants of $X$ to the boundary invariants of the
nilization of $X$.

Next we need some notation on chain complexes $(Y,d)$. Let
$B_nY=dY_{n+1}$ and $Z_nY=\ker\{d:Y_n\to Y_{n-1}\}$ be the modules
of boundaries and cycles respectively. For a chain map $\xi: Y\to
Y'$ we consider sequences of homomorphisms
$$
\delta=(\delta: B_{n-1}Y\to Z_nY')_{n\in \mathbb Z}
$$
and we call the map $\xi+\delta$ obtained by
$(\xi+\delta)_n=\xi_n+i\delta_nd$ the {\it $\delta$-deformation}
of $\xi$, see 4.5.5 \cite{HTH}. The $\delta$-deformation
$1+\delta$ of the identity is an isomorphism of chain complexes
with inverse $1-\delta$. Given an object $(Y,b,\beta)$ in ${\sf
bype}_r(M)$ we get the $\delta$-induced object
$^\delta(Y,b,\beta)=(Y,(1+\delta)_*b,(1+\delta)_*\beta)$ together
with the isomorphism
$$
1+\delta: (Y,b,\beta)\to ^\delta(Y,b,\beta)
$$
in ${\sf bype}_r(M)$. By the following lemma this map $1+\delta$
is always $\lambda$-realizable.
\begin{lemma}\label{le2}
Let $(Y,b,\beta)=\lambda G$ with $G$ a free object in $({\sf
sNil})_r$. Then there is a free object $^\delta G$ in $({\sf
sNil})_r$ together with an isomorphism $^\delta 1:G\to ^\delta G$
such that the composite
$$
(Y,b,\beta)=\lambda G\buildrel{\lambda( ^\delta
1)}\over\rightarrow \lambda( ^\delta G)= ^\delta(Y,b,\beta)
$$
coincides with $1+\delta$.
\end{lemma}
\begin{proof}
We define $^\delta G$ by the following diagram:
$$
\xyma{ & \Lambda^2N^{-1}Y \ar@{->}[r]^{(1+\delta)_*} \ar@{->}[d]
& \Lambda^2N^{-1}Y \ar@{>->}[d] \ar@{=}[r] & \Lambda^2N^{-1}Y \ar@{>->}[d]\\
& G \ar@{->}[r]^\alpha \ar@{->}[d] & G'
\ar@{->}[d] & ^\delta G \ar@{->}[l]^{\beta} \ar@{->}[d]\\
& ab(G) \ar@{=}[r] & N^{-1}Y & N^{-1}Y \ar@{->}[l]^{(1-\delta)_*}\\
}
$$
The columns are short exact and we take the central push out and
the pull back of groups. Then $^\delta 1=\beta^{-1}\alpha$.
\end{proof}
Moreover a $\delta$-deformation has the following property:
\begin{lemma}\label{le3}
Let $Y$ and $Y'$ be chain complexes of free abelian groups and let
$\xi: Y\to Y'$ be a homotopy equivalence and let $\xi+\gamma$ be a
$\gamma$-deformation of $\xi$. Then there is $\delta$ such that
$(1+\delta)\xi$ and $\xi+\gamma$ are homotopic.
\end{lemma}
\begin{proof}
The map $\xi$ induces an isomorphism
$$
\xi_*: Ext(H_nY,H_{n+1}Y)\buildrel{\simeq}\over\rightarrow
Ext(H_nY, H_{n+1}Y'),
$$
where $\gamma_n$ represents an element $\{q\gamma_n\}\in Ext(H_nY,
H_{n+1}Y')$. Let $\delta_n$ be an element which represents
$\{q\delta_n\}=(\xi_*)^{-1}\{q\gamma_n\}\in Ext(H_nY,H_{n+1}Y).$
\end{proof}
Moreover we need for $r\geq 1$ the category ${\sf H}_{n+1}$ of
$r$-reduced homotopy systems of order $n+1$ in 4.2 \cite{HTH}.
Objects in ${\sf H}_{n+1}$ are triple $(C,f_{n+1}, X^n)$ where
$X^n$ is an $n$-dimensional CW-complex with $X^r=*$ and $C$ is a
chain complex of free abelian groups and $f_{n+1}: C_{n+1}\to
\pi_nX^n$ is a homomorphism. These data satisfy the properties in
4.2.2 \cite{HTH}. A morphism in ${\sf H}_{n+1}$ is a pair
$$
(\xi, \eta): X=(C,f_{n+1},X^n)\to Y=(C',g_{n+1},Y^n)
$$
where $\xi: C\to C'$ is a chain map and $\eta: X^n\to Y^n$ is the
0-homotopy class of a cellular map satisfying the properties in
4.2.2 \cite{HTH}. There is a homotopy relation on ${\sf H}_{n+1}$
as in 4.2.6 \cite{HTH} which yields the homotopy category ${\sf
H}_{n+1}/\simeq$. We define in the same way the category ${\sf
H}_{n+1}^{nil}$ by replacing CW-complexes in ${\sf H}_{n+1}$ by
CW-objects (i.e. free objects) in ${\sf sNil}$. Here we are aware
of the fact that $n$-cells in a CW-complex correspond to free
generators of degree $n-1$ in a CW-object in ${\sf sNil}$.\\

\noindent{\bf Remark.} In the book \cite{BauesSpringer} the
category ${\sf H}_{n+1}$ is defined in any cofibration category
with spherical objects. We can apply this to the category ${\sf
sNil}$ since ${\sf sNil}$ has the Blakers-Massey property with
respect to the theory ${\sf T\subset Ho(sNil)}$ given by
coproducts of spherical objects $S(0)$. See Appendix B below. The
tower of categories in \cite{BauesSpringer} shows that ${\sf
H}_{n+1}^{nil}$ has similar properties as $H_{n+1}$ in \cite{HTH},
in particular, the obstructions satisfy formulas as in section 4.5
of \cite{HTH}.
\\

Moreover we define the category ${\sf N}_{n+1}^b$. Objects are
tuple $(C,f_{n+1},X^n,b,\beta)$ where $(C,f_{n+1},X^n)$ is an
object in ${\sf H}_{n+1}^{nil}$ and $(C,b,\beta)$ is an object in
${\sf bype}_r(\mathbb Z^\Lambda)$ such that for $t\leq n+1$ the
homotopy invariants $(b_tX,\beta_{t-1}X)$ of $X=(C,f_{n+1},X^n)$
coincide with $(b_t,\beta_{t-1})$ given by $(b,\beta)$. Morphisms
are morphisms $(\xi,\eta)$ in ${\sf H}_{n+1}^{nil}/\simeq$ for
which $\xi$ is also a morphism in ${\sf bype}_r(\mathbb
Z^\Lambda)$. Compare 4.6.2 \cite{HTH}. The categories ${\sf
N}_{n+1}^b/\simeq$ form a tower of categories, $n\geq r+1$,
$$
{\sf Ho(sNil)}_r \buildrel{r^{n+1}}\over\rightarrow {\sf
N}_{n+1}^b/_\simeq \buildrel{\lambda^{n+1}}\over\rightarrow {\sf
N}_n^b/_\simeq \rightarrow \dots \rightarrow {\sf
N}_{r+1}^b/_\simeq \cong {\sf bype}_r(\mathbb Z^\Lambda)
$$
where the functors $r^{n+1}$ and $\lambda^{n+1}$ are defined as in
4.2.3 \cite{HTH}. The composite of these functors is the functor
$\lambda$ in theorem \ref{mainthe}.
\begin{lemma}
The functor $\lambda:{\sf N}_{n+1}^b/_\simeq \to {\sf
N}_n^b/_\simeq$ is representative. In fact, for an object $X$ in
${\sf N}_n^b$ there is an object $\bar X$ in ${\sf N}_{n+1}^b$
with $\lambda \bar X=X$.
\end{lemma}
\begin{proof}
This is a consequence of 4.4.5 \cite{HTH} where we obtain $\bar X$
in ${\sf H}_{n+1}$. In the same way we obtain $\bar X$ in ${\sf
H}_{n+1}^{nil}$.
\end{proof}
\begin{prop}
Let $r\geq 1$. The functor $\lambda: {\sf Ho(sNil)}_r\to {\sf
bype}_r(\mathbb Z^\Lambda)$ is representative.
\end{prop}
\begin{proof}
Using the lemma we construct inductively for $n\geq r$ objects
$(C,\delta_{n+1},X^n)=X^{(n+1)}$ in ${\sf N}_{n+1}^b$. Then
$\tilde X=\dlim{\strut lim} X^{(n)}$ satisfies
$\lambda(G)=(G,b,\beta)$ with $G$ in ${\sf (sNil)}_r$ satisfying
$\tilde X\simeq B|G|$.
\end{proof}
\begin{lemma}\label{l6}
Let $X$ and $Y$ be objects in ${\sf N}_{n+1}^b$ and let
$F=(\xi,\eta): \lambda X\to \lambda Y$ be a map in ${\sf
N}_n^{b}$. Then there exists $\tau: B_nC\to Z_{n+1}C'$ such that
$(\xi+\tau, \eta): \lambda X\to \lambda Y$ in ${\sf H}_{n}^{nil}$
is $\lambda$-realizable by a map $\bar F: X\to Y$ in ${\sf
H}_{n+1}^{nil}$, that is, $\lambda(\bar F)=(\xi+\tau,\eta)$.
\end{lemma}
\begin{proof}
This is a consequence of 4.6.1 \cite{HTH} transformed to the
category ${\sf H}_{n+1}^{nil}$.
\end{proof}
For an object $Y=(C',g_{n+1},Y^n,b',\beta')$ in ${\sf N}_{n+1}^b$
and for $\delta: B_nC'\to Z_{n+1}C'$ let $^\delta
Y=(C',g_{n+1},Y^n, (1+\delta)_*b',(1+\delta)_*\beta').$
\begin{lemma}\label{le7}
Let $X$ and $Y$ be objects in ${\sf N}_{n+1}^b$ and let
$F=(\xi,\eta): \lambda X\to \lambda Y$ be a map in ${\sf N}_n^b$
where $\xi$ is a homotopy equivalence of chain complexes. Then
there exists $\delta: B_nC'\to Z_{n+1}C'$ such that
$$
((1+\delta)\xi,\eta):\lambda X\to \lambda(^\delta Y)
$$
in ${\sf N}_n^b/_\simeq$ is $\lambda$-realizable by a map $\bar
F:X\to ^\delta Y$ in ${\sf N}_{n+1}^b/_\simeq$, that is
$\lambda(\bar F)=((1+\delta)\xi,\eta)$.
\end{lemma}
\begin{proof}
We use lemma \ref{le3} and lemma \ref{l6}.
\end{proof}
\begin{prop}
Let $r\geq 1$. The restriction of $\lambda: {\sf Ho(sNil)}_r\to
{\sf bype}_r(\mathbb Z^\Lambda)$ to the subcategories of
isomorphisms is a full functor.
\end{prop}
\begin{proof}
Let $G$ and $G'$ be free objects in ${\sf (sNil)}_r$ with
$C=N^{-1}ab(G)$ and $C'=N^{-1}ab(G')$. Let
$$
\xi: (C,b,\beta)=\lambda G\to (C',b',\beta')=\lambda G'
$$
be a homotopy equivalence in ${\sf bype}_r(\mathbb Z^\Lambda)$.
Let $X(n)$ be the objects in ${\sf N}_n^b/_\simeq$ given by $G$
and the tower of categories above and by $b,\beta$ in
$\lambda(G)$. Then $\xi$ determines a map $\xi: X(r+1)\to Y(r+1)$
in ${\sf N}_{r+1}^b$ which by lemma \ref{le7} yields a map
$$
\bar F: X(r+2)\to ^{\delta_{r+1}}Y(r+2)
$$
in ${\sf N}_{r+2}^b$ which again by lemma \ref{le7} yields a map
$$
\bar{\bar F}: X(r+3)\to ^{\delta_{r+2}}(^{\delta_{r+1}}Y(r+3))
$$
in ${\sf N}_{r+3}^b$ and so on. Inductively we get a map in ${\sf
sNil}$
$$
F^\infty: X=G\to ^\delta Y= ^\delta G'
$$
Here $F^\infty$ determines a map $F:G\to ^\delta G'$ in ${\sf
Ho(sNil)}$ for which $\lambda(F)=(1+\delta)\xi$. Now lemma
\ref{le2} shows that there is a map $F':G\to G'$ with $\lambda
F'=\xi$.
\end{proof}
 \vspace{.5cm}
\section{Proof of Theorem \ref{lamq}}

\vspace{.5cm} The formulas (I)-(IV) and (XI) directly follow from
the universal coefficient theorem for the exterior square functor
(\ref{deri}) together with (\ref{byper}). For the computation of
other functors, recall the following construction due to Curtis
\cite{SH}. Let $X$ be a simplicial group. The lower central series
filtration in $X$ gives rise to the long exact sequence
\begin{multline*}
\dots\to \pi_{i+1}(X/\gamma_n(X))\to
\pi_i(\gamma_n(X)/\gamma_{n+1}(X))\to\\
\pi_i(X/\gamma_{n+1}(X))\to \pi_i(X/\gamma_n(X))\to \dots
\end{multline*}

This exact sequence defines the graded exact couple which gives
rise to the natural spectral sequence $E(X)$ with the initial
terms
\begin{align*}
& E_{p,q}^1(X)=\pi_q(\gamma_{p}(X)/\gamma_{p+1}(X)).
\end{align*}
and the differentials $d^i,\ i\geq i$
\begin{align*}
& d^i: E_{p,q}^i(X)\to E_{p+i, q-1}^i(X).
\end{align*}
\begin{lemma}\cite{SH}\label{curtishom}
Let $K$ be a connected and simply connected simplicial set, $G=GK$
its Kan's construction. Then the spectral sequence $E^i(G)$
converges to $E^\infty(G)$ and $\oplus_rE_{p,q}^\infty$ is the
graded group associated with the filtration on
$\pi_q(GK)=\pi_{q+1}(|K|)$. The groups $E^1(K)$ are homology
invariants of $K$.
\end{lemma}
For a given abelian group $D$, and Eilenberg-MacLane space
$K(D,m),\ m\geq 1$, consider the Kan construction $GK(D,m)$ and
corresponding Curtis spectral sequence
\begin{equation}\label{conve}
E_{p,q}^1(GK(D,m))\Rightarrow \pi_{q+1}K(D,m)=\begin{cases} D,\
q=m-1,\\
0,\ \text{otherwise}\end{cases}.
\end{equation}
Denote $Y=ab(GK(D,m))$. We naturally have
$$
E_{1,m-1}^1(GK(D,m))=E_{1,m-1}^\infty(GK(D,m))=\pi_{m-1}(Y)=H_mK(D,m)=D.
$$
The convergence (\ref{conve}) means that
$$
E_{p,q}^\infty=0,\ (p,q)\neq (1,m-1).
$$
Recall also the connectivity result due to Curtis
\cite{Curtis:63}. Let $X$ be a free simplicial group, which is
$m$-connected ($m\geq 0$). Then $\gamma_r(X)/\gamma_{r+1}(X)$ is
$\{m+log_2r\}$-connected, where $\{a\}$ is the least $\geq a$.
Applying this result to our case, we get that
$$\pi_i(L^rY)=0,\ i<m+2,\ r\geq 3$$
where $L^r: {\sf Ab\to Ab}$ is the $r$-th Lie functor. Collect the
low dimensional elements of the spectral sequence
$E_{p,q}^1(GK(D,m))$ in
the following table:\\
\begin{center}
\begin{tabular}{ccccccc}
$m+2$ & \vline & $H_{m+3}K(D,m)$ & $\pi_{m+2}(\Lambda^2Y)$ & $\pi_{m+2}(L^3Y)$ & $\pi_{m+2}(L^4Y)$ & $\pi_{m+2}(L^5Y)$\\
$m+1$ & \vline & $H_{m+2}K(D,m)$ & $\pi_{m+1}(\Lambda^2Y)$ & 0 & 0 & 0\\
$m$ & \vline & $H_{m+1}K(D,m)$ & $\pi_{m}(\Lambda^2Y)$ & 0 & 0 & 0\\
$m-1$ & \vline & $D$ & $\pi_{m-1}(\Lambda^2Y)$ & 0 & 0 & 0\\
\hline
$q/ p$ & \vline & 1 & 2 & 3 & 4 & 5\\
\end{tabular}
\end{center}
\vspace{.5cm} The convergence (\ref{conve}) then implies that the
differentials
\begin{equation}\label{isor}
d_{1,i}^1: H_{i+1}K(D,m)\to \pi_{i-1}(\Lambda^2Y),\ i=m,m+1,m+2
\end{equation}
are isomorphisms.

The simplicial fibration sequence
$$
\Lambda^2 Y\to nil(GK(D,m))\to Y
$$
induces the map of chain complexes
$$
f: NY[1]\to N\Lambda^2 Y,
$$
which induces the isomorphisms of homology groups (\ref{isor}) in
dimensions $m, m+1, m+2$. Hence for every abelian group $A$, one
has natural isomorphism
$$
[C(A,i),NY]\simeq [C(A,i-1),N\Lambda^2Y], i=m,m+1.
$$
Since $H_{m+1}K(D,m)=0$, the natural map
$$
Y\to N^{-1}C(D,m-1)
$$
induces isomorphisms of homology groups
$$
H_iN\Lambda^2Y\to H_i\mathbb Z_\#^\Lambda(D,m-1)
$$
for $i\leq m+2$ and therefore
$$
[C(A,i),NY]\simeq [C(A,i-1),\mathbb Z_\#^\Lambda C(D,m-1)],\
i=m,m+1.
$$
Recall now the definition of the pseudo-homology functors
$H_n^{(m)}$ (see 6.3 \cite{HTH}). Given abelian groups $A,D$,
$$
H_n^{(m)}(A,D)=H_n(A,K(D,m)).
$$
In the simplicial language we have
$$
H_n^{(m)}(A,D)=[C(A,n-1),Nab(GK(D,m))].
$$
Therefore
$$
H_{m+2}^{(m)}(A,D)=[C(A,m+1), Nab(GK(D,m))]=[C(A,m),\mathbb
Z_\#^\Lambda C(D,m-1)]=Sq_{m,m-1}(A,D).
$$
There is a certain periodicity in the description of homotopy type
of the simplicial group $\Lambda^2N^{-1}C(D,m)$ reflected in the
formulation of the universal coefficient theorem for quadratic
modules. This principle directly implies that
\begin{align*}
& Sq_{2m,m}^{\mathbb Z^\Lambda}(A,D)=H_4^{(2)}(A,D),\ m\ \text{odd},\\
& Sq_{2m-1,m}^{\mathbb Z^\Lambda}(A,D)=H_5^{(3)}(A,D),\ m\ \text{even},\\
& Sq_{n,m}^{\mathbb Z^\Lambda}(A,D)=H_6^{(4)}(A,D),\ n-m\
\text{odd},\ m<n<2m-1.
\end{align*}
For the pseudo-homology functors $H_{n}^{(m)}$ one has the direct
sum decompositions (see 6.3.9 \cite{HTH}). In the cases which we
need here the decompositions are the following:
\begin{align*}
& H_4^{(2)}(A,D)=\Gamma T_\#(A,D)\\
& H_5^{(3)}(A,D)=L_\#(A,D)\oplus Ext(A, D*\mathbb Z_2)\\
& H_6^{(4)}(A,D)=Ext(A, D*\mathbb Z_2)\oplus Hom(A,D\otimes\mathbb
Z_2)
\end{align*}
Therefore the formulas (VI), (VIII) and (X) follow.

For the description of other functors $Sq_{n,m}^{\mathbb
Z^\Lambda}$, consider the following diagram with exact rows:
\begin{equation}\label{diagra1} \xyma{ Ext(A, H_{m+4}K(D,m))
\ar@{>->}[r] \ar@{->>}[d]^{(d_{1,m+3}^1)^*} & H_{m+3}^{(m)}(A,D)
\ar@{->}[d]
\ar@{->>}[r] & Hom(A, H_{m+3}K(D,m)) \ar@{->}[d]^{(d_{1,m+2}^1)^*} \\
Ext(A,\pi_{m+2}\Lambda^2Y) \ar@{>->}[r] \ar@[->>][d] &
[C(A,m+1),\Lambda^2Y] \ar@{->>}[r] \ar@{->}[d] &
Hom(A,\pi_{m+1}\Lambda^2Y) \ar@{->}[d]\\ Ext(A,\pi_{m+2}\mathbb
Z_\#^\Lambda C(D,m-1)) \ar@{>->}[r] & Sq_{m+1,m-1}(A,D)
\ar@{->>}[r] & Hom(A,\pi_{m+1}\mathbb Z_\#C(D,m-1))
\\ }
\end{equation}
The above argument shows that the right hand vertical arrows in
(\ref{diagra1}) are isomorphisms. Hence, the functors
$Sq_{n,m}^{\mathbb Z^\Lambda}$ for certain $m,n$ in the diagram
(\ref{diagra1}) can be described as the following push-outs:
\begin{align*}
& \xyma{ & Ext(A, D\otimes (\mathbb Z_2\oplus \mathbb Z_3))
\ar@{>->}[r]
\ar@{->>}[d] & H_8^{(5)}(A,D) \ar@{->}[d] \ar@{->>}[r] & Hom(A, D*\mathbb Z_2)\\
& Ext(A,D\otimes \mathbb Z_2) \ar@{->}[r] & Sq_{6,4}^{\mathbb Z^\Lambda}(A,D) & \\
}\\ &  \xyma{ & Ext(A, \Gamma(D)\oplus D\otimes \mathbb Z_3)
\ar@{>->}[r]
\ar@{->>}[d] & H_7^{(4)}(A,D) \ar@{->}[d] \ar@{->>}[r] & Hom(A, D*\mathbb Z_2)\\
& Ext(A,\Gamma(D)) \ar@{->}[r] & Sq_{5,3}^{\mathbb Z^\Lambda}(A,D) & \\
}\\ & \xyma{ & Ext(A, \Omega(D)\oplus D\otimes \mathbb Z_3)
\ar@{>->}[r]
\ar@{->>}[d] & H_6^{(3)}(A,D) \ar@{->}[d] \ar@{->>}[r] & Hom(A, \Lambda^2(D)\oplus D*\mathbb Z_2)\\
& Ext(A,\Omega(D)) \ar@{->}[r] & Sq_{4,2}^{\mathbb Z^\Lambda}(A,D) & \\
}
\end{align*}
The decompositions for functors $H_{m+3}^{(m)}$ (see 6.3.10
\cite{HTH}) and the periodicity principle imply the following
decompositions of the needed functors:
\begin{align*}
& Sq_{n,m}^{\mathbb Z^\Lambda}(A,D)=Ext(A, D\otimes \mathbb Z_2)\oplus Hom(A,D*\mathbb Z_2),\ n-m\ \text{even},\ m<n<2n-1,\\
& Sq_{2m-1,m}^{\mathbb Z^\Lambda}(A,D)=Ext(A, \Gamma(D))\oplus Hom(A,D*\mathbb Z_2),\ m\ \text{odd},\\
& Sq_{2m,m}^{\mathbb Z^\Lambda}(A,D)=\Lambda^2T_\#(A,D)\oplus
Hom(A, D*\mathbb Z_2),\ m\ \text{even}
\end{align*}
and the formulas (V), (VII) and (IX) follow.\ $\Box$
 \vspace{.5cm}
\section*{Appendix A: Homotopy groups of spherical objects in categories of nilpotent groups}
\vspace{.5cm} The homotopy theory in the category $\sf sNil$ in
this paper can be generalized for higher nilpotency degree. In
particular, it is possible to compute homotopy groups of spherical
objects in the category of simplicial $r$-nilpotent groups for
$r=3$ and partially for $r=4,5$. We give examples as follows.

Let $n\geq 1,\ r\geq 2$. Consider the category ${\sf Nil}^r$ of
nilpotent groups of class $r$. The homotopy groups of $S^{n+1}$ in
the category ${\sf Nil}^r$ can be naturally defined via Milnor's
$F[S^n]$-construction:
$$
\pi_{i+1}^{{\sf Nil}^r}(S^{n+1})=\pi_iF[S^n]_{{\sf
Nil}^r}=\pi_i(F[S^n]/\gamma_{r+1}(F[S^n])).
$$
Here $G_{{\sf Nil}^r}$ is the $r$-nilization of a simplicial group
$G$ in the category ${\sf Nil}^r$. Denote $K(\mathbb
Z,n)=ab(F[S^n])$. Here we consider low-dimensional cases. As we
will see, the general problem of the description of homotopy
groups in categories ${\sf Nil}^r$ essentially reduces to the
homotopical properties of the $r$-th
Lie functor $L^r:\ {\sf Ab}\to {\sf Ab}.$\\

\noindent {\it Category ${\sf Nil^2}$}. The homotopy groups of
spheres in ${\sf Nil}^2={\sf Nil}$ are computed in $\S 4$. Recall
them:
\begin{equation}\label{nil2sphere}
\pi_{n+k+1}^{{\sf Nil}^2}(S^{n+1})=\begin{cases} \mathbb Z,\ k=0\\
\mathbb Z,\ k=n\ \text{odd}\\
\mathbb Z_2,\ 0<k<n,\ k\ \text{odd}\\
0,\ \text{otherwise}\end{cases}
\end{equation}

\noindent{\it Category ${\sf Nil^3}$}: There is the following
natural exact sequence of simplicial groups
$$
1\to L^3K(\mathbb Z,n)\to F[S^n]_{{\sf Nil}^3}\to F[S^n]_{{\sf
Nil}^2}\to 1
$$
which induces the long exact sequence of corresponding homotopy
groups. We will use the following result due to Schlesinger
\cite{Sch}: if $p$ is an odd prime then
\begin{equation}\label{primedec}
\pi_{n+k}L^pK(\mathbb Z,n)=\begin{cases} \mathbb Z_p,\
k=2i(p-1)-1,\ i=1,2,\dots, [n/2]\\
0, \text{otherwise}\end{cases}
\end{equation}
Hence
\begin{equation}\label{dfr}
\pi_{n+k}L^3K(\mathbb Z,n)=\begin{cases} \mathbb Z_3,\ k=4i-1,\
i=1,2,\dots, [n/2]\\
0,\ \text{otherwise}
\end{cases}
\end{equation}
The description of the homotopy groups in category ${\sf Nil}^2$
(\ref{nil2sphere}) and (\ref{dfr}) then imply the following:
\begin{equation}\label{nil3sphere}
\pi_{n+k+1}^{{\sf Nil^3}}(S^{n+1})=\begin{cases} \mathbb Z,\ k=0,\\
\mathbb Z_2,\ 0<k<n,\ k\equiv 1\mod 4,\\
\mathbb Z_6,\ 0<k<n,\ k\equiv 3\mod 4,\\
\mathbb Z_3,\ n<k<2n,\ k\equiv 3\mod 4\\
\mathbb Z\oplus \mathbb Z_3,\ n\equiv 3\mod 4,\ k=n\\
\mathbb Z,\ n\equiv 1\mod 4,\ k=n\\
\end{cases}
\end{equation}

 \noindent {\it Category ${\sf Nil^4}$}: We have the short exact sequence of
 simplicial groups
$$
1\to L^4K(\mathbb Z,n)\to F[S^n]_{{\sf Nil}^4}\to F[S^n]_{{\sf
Nil}^3}\to 1
$$
The description (\ref{nil3sphere}) essentially reduces the problem
of computation of $\pi_*^{\sf Nil^4}(S^n)$ to homotopical
properties of simplicial abelian groups $L^4K(\mathbb Z,n)$. There
is the following natural decomposition of the fourth Lie functor
$L^4$ for a free $\mathbb Z$-module $M$:
$$
0\to \Lambda^2\Lambda^2(M)\to L^4(M)\to J^4(M)\to 0,
$$
where $J^4$ is the fourth metabelian Lie functor, which can be
defined as the kernel of the symmetrization map:
\begin{equation}\label{symmmap}
0\to J^4(M)\to M\otimes SP^3(M)\to SP^4(M)\to 0
\end{equation}
(here $SP^*$ is the symmetric tensor power). This is a simplest
case of the Curtis decomposition of Lie functors (see
\cite{Curtis:63}).

Homotopy groups of simplicial abelian groups
$\Lambda^2\Lambda^2K(\mathbb Z,n)$ follow from the split exact
sequence (\ref{deri}) (we use the sequence (\ref{deri}) twice). At
the first step we have
$$
\pi_i\Lambda^2\Lambda^2K(\mathbb Z,1)=\begin{cases} \mathbb
 Z_2,\ i=3\\
 0,\ \text{otherwise}\end{cases}
$$
For the 2-sphere one has a contractible $J^4K(\mathbb Z,1)$ (see
\cite{Curtis:63}), therefore there is a weak homotopy equivalence
$$
\Lambda^2\Lambda^2K(\mathbb Z,1)\buildrel{\sim}\over\rightarrow
L^4K(\mathbb Z,1)
$$
and the following description of the homotopy groups follows
immediately:
\begin{equation}\label{sph2}
\pi_i^{{\sf Nil}^4}(S^2)=\begin{cases} \mathbb Z,\ i=2,3\\
\mathbb Z_2,\ i=4,\\
0,\ \text{otherwise}
\end{cases}
\end{equation}
The case of the 3-sphere is more complicated. First of all we have
$$
\pi_i\Lambda^2\Lambda^2K(\mathbb Z,2)=\begin{cases} \mathbb Z_4,\
i=6\\
\mathbb Z_2,\ i=4,5,7\\
0,\ \text{otherwise}
\end{cases}
$$
However in this case the functor $J^4$ also contributes. It
directly follows from \cite{DT} and \cite{DP} (see p.307) that
$J^4K(\mathbb Z,2)=K(\mathbb Z_4,7).$ The short exact sequence of
functors (\ref{symmmap}) implies the following boundary
homomorphism
$$
\begin{CD}
\pi_7J^4K(\mathbb Z,2) @>\partial>>
\pi_6\Lambda^2\Lambda^2K(\mathbb Z,2)\\
@| @|\\
\mathbb Z_4 @>>> \mathbb Z_4
\end{CD}
$$
with $\pi_6L^4K(\mathbb Z,2)=coker(\partial)$. The direct
simplicial computations show that $\partial$ is an isomorphism. As
a result we obtain the following
\begin{equation}\label{sph3}
\pi_i^{{\sf Nil}^4}(S^3)=\begin{cases} \mathbb Z,\ i=3\\
\mathbb Z_2,\ i=4,5,8\\
\mathbb Z_6,\ i=6\\
0,\ \text{otherwise}
\end{cases}
\end{equation}
The general description of the homotopy groups $\pi_*^{\sf
Nil^4}(S^n)$ seems to be quite nontrivial and we leave it as an
open problem. Observe that (\ref{sph2}), (\ref{sph3}) together
with (\ref{primedec}) imply that
$$
\pi_i^{{\sf Nil}^5}(S^2)=\begin{cases} \mathbb Z,\ i=2,3\\
\mathbb Z_2,\ i=4,\\
0,\ \text{otherwise}
\end{cases}
$$
and
$$
\pi_i^{{\sf Nil}^5}(S^3)=\begin{cases} \mathbb Z,\ i=3\\
\mathbb Z_2,\ i=4,5,8\\
\mathbb Z_6,\ i=6\\
\mathbb Z_5,\ i=10\\ 0,\ \text{otherwise}
\end{cases}
$$
\vspace{.5cm}
\section*{Appendix B: Blakers-Massey property for {\sf sNil}}\vspace{.5cm} Let $\sf C$ be
the model category $\sf sGr$ or ${\sf sNil}^k$ for $k\geq 1$ and
let $\sf T\subset {\sf Ho(C)}$ be the theory of coproducts of
$0$-spheres in $\sf C$. According to \cite{BauesSpringer} the
category $\sf C$ under $\sf T$ has the {\it Blakers-Massey
property} if the following condition is satisfied. Let $F, M, N$
be cofibrant objects in ${\sf C}$ such that $M$ and $N$ are
obtained from $F$ by attaching cells in dimensions $\geq m$ and
$\geq n$ respectively, in particular there are monomorphisms $f_1:
F\to M$ and $f_2: F\to N$, which preserve bases. Consider the push
out in ${\sf C}$:
\begin{equation}\label{pushka}
\begin{CD}
N @>>> M\vee_FN\\
@A{f_2}AA @AAA\\
F @>{f_1}>> M
\end{CD}
\end{equation}
For ${\sf C=sNil}^k$ the push out $M\vee_F N$ can be obtained as a
nilization of the amalgamated product $M*_{F}N$. Then the induced
map of relative homotopy groups
\begin{equation}\label{abpush}
g_r: \pi_r(N, F)\to \pi_r(M\vee_F N, M)
\end{equation}
is an isomorphism for $r\leq n+m-2$ and an epimorphism for
$r=n+m-1$. The classical Blakers-Massey theorem implies that the
Blakers-Massey property holds for $\sf C=sGr$.

\begin{theorem}\label{bmp}
For ${\sf C=sNil}^k,\ k\geq 1$, the Blakers-Massey property holds.
\end{theorem}

The theorem implies that all the theory of \cite{BauesSpringer}
can be applied to the category ${\sf sNil}^k$. In particular one
gets
\begin{cor}
There is the tower of categories as in \cite{BauesSpringer} for
${\sf C=sNil}^k$.
\end{cor}

\noindent{\it Proof of Theorem \ref{bmp}.} We prove the theorem
for $k=2$, the general proof can be obtain from the given one
using the inductive arguments. Suppose first that the push out
(\ref{pushka}) is considered in the category ${\sf sAb}$ of
simplicial abelian groups. Since one has a natural isomorphism of
abelian groups
$$
\pi_i(N,F)\simeq H_i(N/F),\ \pi_i(M\vee_{F}N,M)=H_i(N/F),
$$
the push-out (\ref{pushka}) induces natural isomorphisms
(\ref{abpush}) for all $i>0$. Now consider (\ref{pushka}) in the
category ${\sf sNil^2}$. Since the maps $f_1,f_2$ preserve bases,
the diagram (\ref{pushka}) induces an analogous diagram of
abelianizations and we have the following commutative diagrams:
\begin{equation}\label{j4}
\begin{CD}
\pi_i(\Lambda^2abF) @>>> \pi_i(\Lambda^2abN) @>>>
\pi_i(\Lambda^2abN/\Lambda^2abF) @>>> \pi_{i-1}(\Lambda^2abF)\\
@VVV @VVV @VVV @VVV\\
\pi_i(F) @>>> \pi_i(N) @>>> \pi_i(N,F) @>>> \pi_{i-1}(F)\\
@VVV @VVV @VVV @VVV\\
\pi_i(abF) @>>> \pi_i(abN) @>>> \pi_i(abN/abF) @>>> \pi_{i-1}(abF)
\end{CD}
\end{equation}
and
\begin{equation}\label{j5}
\xyma{\pi_i(\Lambda^2abM)) \ar@{->}[r] \ar@{->}[d] &
\pi_i(\Lambda^2ab(M\vee_FN)) \ar@{->}[r] \ar@{->}[d] &
\pi_i(\Lambda^2ab(M\vee_FN)/\Lambda^2abM) \ar@{->}[r] \ar@{->}[d]
& \pi_{i-1}(\Lambda^2abM) \ar@{->}[d]\\
\pi_i(M) \ar@{->}[r] \ar@{->}[d] & \pi_i(M\vee_FN) \ar@{->}[r]
\ar@{->}[d] & \pi_i(M\vee_FN,M) \ar@{->}[r] \ar@{->}[d] &
\pi_{i-1}(M) \ar@{->}[d] \\
\pi_i(abM) \ar@{->}[r] & \pi_i(ab(M\vee_FN)) \ar@{->}[r] &
\pi_i(abN/abF) \ar@{->}[r] & \pi_{i-1}(abM)\\ }
\end{equation}
A general argument in model categories implies that all rows and
columns in (\ref{j4}) and (\ref{j5}) are exact.  The above
diagrams induce the natural commutative diagram
\begin{equation}\label{oka1}
\begin{CD}
\pi_i(\Lambda^2abN/\Lambda^2abF) @>{j_i}>>
\pi_i(\Lambda^2ab(M\vee_FN)/\Lambda^2abM)\\
@VVV @VVV\\
\pi_i(N,F) @>{g_i}>> \pi_i(M\vee_FN,M)\\
@VVV @VVV\\
\pi_i(abN/abF) @= \pi_i(abN/abF)
\end{CD}
\end{equation} Furthermore we have the following commutative diagrams
\begin{equation}\label{oka2}
\begin{CD}
\pi_i(abF\otimes (abN/abF)) @>{p_i}>> \pi_i(abM\otimes
(ab(M\vee_FN)/abM))\\
@VVV @VVV\\ \pi_i(\Lambda^2abN/\Lambda^2abF) @>{j_i}>>
\pi_i(\Lambda^2ab(M\vee_FN)/\Lambda^2abM)\\
@VVV @VVV\\ \pi_i(\Lambda^2(abN/abF)) @=
\pi_i(\Lambda^2(ab(M\vee_FN)/abM))
\end{CD}
\end{equation}
where the maps
$$
p_i: \pi_i(abF\otimes (abN/abF))\to \pi_i(abM\otimes
(ab(M\vee_FN)/abM))=\pi_i(abM\otimes (abN/abF))
$$
are induced by inclusion $f_1$. Since $N$ is obtained from $F$ by
adding elements in dimensions $\geq n$, $\pi_i(abN/abF)=0,\ i<n$.
By the same argument, the map $\pi_i(abF)\to \pi_i(abM)$ is an
isomorphism for $i<m-1$ and an epimorphism for $i=m-1$. Therefore,
$p_i$ is an isomorphism for $i\leq m+n-2$ and an epimorphism for
$i=m+n-1$. Diagrams (\ref{oka1}) and (\ref{oka2}) imply that the
same property holds for $j_i$ and $g_i$.\ $\Box$

\end{document}